%% file: 0_main.tex
\documentclass[preprint,12pt]{elsarticle}

\usepackage{amssymb}
\usepackage{amsmath}
\usepackage{amsthm}

\usepackage{tikz}
\usepackage{subcaption}
\usepackage{pgfplots}
\usepackage{pgfplotstable}
\pgfplotsset{compat=1.18}
\usepackage{lmodern}
\usepackage{siunitx}
\usetikzlibrary{calc, backgrounds} 

\usepackage{booktabs}
\usepackage{rotating}
\usepackage{pdflscape}

\usepackage{threeparttable}

\usepackage{url}

\usepackage{longtable}
\usepackage{caption}
\usepackage{csvsimple}

\newdefinition{definition}{Definition}

\journal{Computers \& Operations Research}

\begin{document}

\begin{frontmatter}


\title{
Modified Dynamic Programming Algorithms \\
for Order Picking in Single-Block \\
and Two-Block Rectangular Warehouses
} 


\affiliation[label1]{
    organization={
    School of Information and Physical Sciences, 
    University of Newcastle
    },
    state={NSW},
    country={Australia}
}
\affiliation[label2]{
    organization={
    Department of Industrial and Management Systems Engineering,\\
    University of South Florida
    },
    state={FL},
    country={USA}
}
\affiliation[label3]{
    organization={
    Carey Business School,
    Johns Hopkins University
    },
    state={MD},
    country={USA}
}
\affiliation[label4]{
    organization={
    International Computer Science Institute,\\
    University of California at Berkeley
    },
    state={CA},
    country={USA}
}

\author[label1]{George~Dunn\corref{cor1}}
\ead{george.dunn@uon.edu.au}

\author[label1]{Elizabeth~Stojanovski}

\author[label1]{Bishnu~Lamichhane}

\author[label2]{Hadi~Charkhgard}

\author[label3,label4]{Ali~Eshragh}

\cortext[cor1]{Corresponding author}


\begin{abstract}
Recent research has shown that optimal picker tours
in rectangular warehouses exhibit deterministic travel patterns
within each aisle,
and that certain previously considered traversals
are unnecessary.
Using these insights,
this paper proposes modifications to dynamic programming
algorithms that improve computational
efficiency without affecting optimality.
For layouts with and without a central cross-aisle,
the modifications preserve linear-time complexity
in the number of aisles while reducing the number
of state–action evaluations per stage.
The proposed modifications reduce computational
effort by factors up to 1.81,
confirmed by numerical experiments.
These findings are encouraging and highlight
how structural refinements can yield significant
improvements in practical performance of algorithms.
\end{abstract}







\end{frontmatter}



\input{1_introduction}


\input{2_background}


\input{3_single}



\input{4_two}



\input{5_discussion}




\newpage

\section*{Acknowledgments}

George was supported by an Australian Government
Research Training Program (RTP) Scholarship.

\section*{Data and Code Availability}

All code required to reproduce the computational results
presented in this paper is available on request.


\bibliographystyle{elsarticle-num-names} 
\bibliography{bibliography}


\clearpage

\begin{landscape}

\appendix
\section{Two-block state transitions}
\label{app1}

\input{table_two_v2}

\end{landscape}


\end{document}

%% file: 1_introduction.tex

\section{Introduction}
\label{introduction}


The order picking problem aims to determine
the shortest possible route through a warehouse
that visits all necessary item locations.
The standard warehouse layout has parallel aisles,
a central depot, and
cross-aisles to allow for travel between aisles.
\citet{ratliff1983order} presented a
dynamic programming algorithm
for finding the minimal route for warehouses
with a front and back cross-aisle (single-block),
which was extended by \citet{roodbergen2001routing}
for layouts with an additional central cross-aisle
(two-block).
Both methods are based on sequentially building
a tour by applying a limited set of
valid travel options within each aisle and cross-aisle.


A warehouse layout is said to be rectangular
if all aisles are parallel,
all cross-aisles are parallel,
aisles are perpendicular to cross-aisles,
and the distance between any two consecutive aisles
and any two consecutive cross-aisles is constant.
For such layouts,
recent work has substantially clarified the structure
of optimal tours.
In particular,
prior results have shown that double traversals
of aisles are unnecessary in single-block
and two-block rectangular warehouses and,
more generally,
that connectivity does not require double traversals
in warehouses with any number of cross-aisles
\cite{revenant2025note, dunn2025strict, dunn2025double}.
Furthermore,
it has been established that the optimal
travel within an aisle is uniquely determined
by cross-aisle configurations \cite{dunn2025deterministic}.


This paper builds directly on these results
and shows how they can be exploited algorithmically
to reduce the number of stages and transitions required
by existing dynamic programming formulations.
Specifically:
\begin{itemize}
    \item For rectangular single-block warehouses,
    we propose a modified dynamic programming algorithm that
    preserves linear-time complexity while reducing the number
    of state–action evaluations,
    achieving a theoretical reduction of approximately 1.81
    and matching empirical speedups.
    \item For rectangular two-block warehouses,
    we introduce two successive modifications to the
    dynamic programming algorithm that maintain linear-time
    complexity and reduce computational effort,
    with predicted reductions of 1.18 and 1.31,
    confirmed by computational experiments.
\end{itemize}


The remainder of the paper is structured as follows.
Section \ref{background} provides background on warehouse layouts,
graph representations, and prior results.
Sections \ref{single} and \ref{two}
present the single-block and two-block
dynamic programming algorithms,
including the proposed modifications,
complexity analysis,
and computational experiments.
Section \ref{discussion} concludes the paper with
a summary of contributions
and directions for future research.


%% file: 2_background.tex

\section{Background}
\label{background}


We consider rectangular warehouses with
$m \geq 1$ vertical pick aisles
and $n = 2$ (single-block) or
$n = 3$ (two-block)
horizontal cross-aisles
as illustrated in Figure \ref{fig:warehouse}.
In the two-block layout,
the middle cross-aisle divides each aisle
into lower and upper subaisles.
The aisles are considered narrow enough
so that the horizontal distance needed to move
from one side to the other can be regarded negligible.
The single-block (respectively, two-block)
warehouse can be modeled as a graph $G$.
For each aisle $j \in [1,m]$,
the vertex set includes intersections with the cross-aisles.
In the single-block layout,
vertices $a_j$ and $b_j$ correspond to the intersections
of aisle $j$ with the back and front cross-aisles,
respectively.
In the two-block layout,
vertices $a_j$, $b_j$ and $c_j$
correspond to the intersections of aisle
j with the back, middle, and front cross-aisles,
respectively.
The set of vertices
$P = \{p_0, p_1, ..., p_k\}$
denotes the locations to be visited,
where $p_0$ is the depot and
$p_1, ..., p_k$ correspond to the location of
products that must be collected.
Only $p_0$ may lie on a cross-aisle vertex,
whereas all remaining vertices of $P$
must be situated inside the subaisles.

\input{figure_warehouse}


A subgraph $T \subset G$ is called a tour subgraph
if it includes all vertices of $P$
and there is a picking tour that traverses
every edge of $T$ exactly once.
Thus, determining a solution to the picker routing problem
is equivalent to identifying a tour subgraph
whose total edge length is minimal.
Traditional algorithms achieve this by building
partial tour subgraphs (PTSs),
which are subgraphs for which there exists
at least one completion that forms a valid tour subgraph.
Proceeding sequentially from left to right,
PTSs are generated using a restricted set
of valid vertical edge configurations within each aisle
and horizontal edge configurations between adjacent aisles.


\citet{ratliff1983order} demonstrated that,
for a minimal tour in a single-block warehouse,
there exist six distinct vertical edge configurations,
which are illustrated in Figure \ref{fig:vertical}.
These consist of
$(i)$ a single traversal of the aisle ($1pass$),
$(ii)$ entering and exiting via the back cross-aisle ($top$),
$(iii)$ entering and exiting via the front ($bottom$), 
$(iv)$ traversing the aisle in a way that leaves the largest
section untraveled ($gap$),
$(v)$ traversing the aisle twice
($2pass$, also referred to as a double traversal), and
$(vi)$ not entering the aisle at all ($none$).
Note that $none$ is only valid if there are no
items to be collected within an aisle.
These six configurations extend to
two-block warehouses,
where they apply independently to each subaisle
\cite{roodbergen2001routing}.

\input{figure_vertical}


Subsequent research has shown that
the set of valid vertical configurations
considered in earlier algorithms can be reduced.
\citet{revenant2025note} proved that $2pass$
is not required in optimal tours for
rectangular single-block warehouses.
This result was recently extended by the present authors
to two-block warehouses with more than one non-empty aisle
\cite{dunn2025strict}.
More generally,
\citet{dunn2025double} showed that double traversals
are not required to ensure connectivity of the tour
in rectangular warehouses with any number of cross-aisles.
Based on this structural result,
\citet{dunn2025deterministic} established that
once the horizontal edge configuration of a warehouse is fixed,
the corresponding vertical edge configuration
is uniquely determined.
These results form the theoretical foundation
for the algorithmic modifications proposed
in the remainder of this paper.

%% file: figure_warehouse.tex
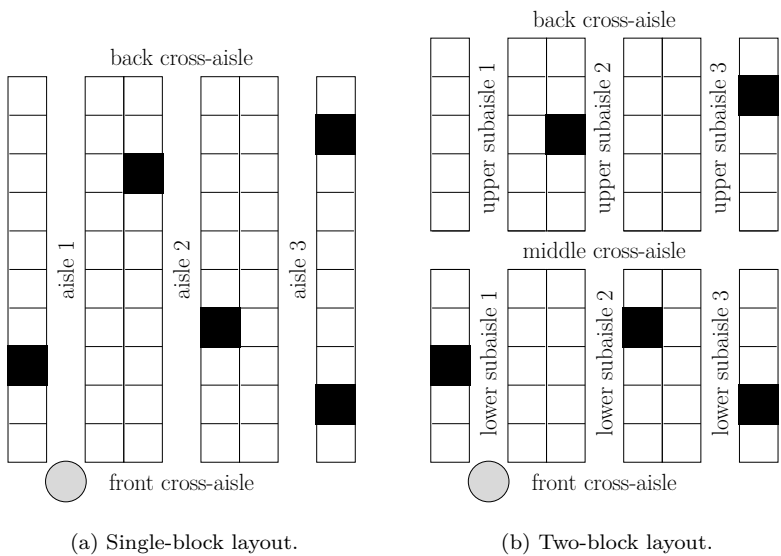
\begin{figure}[ht]
\vskip 0.2in
\begin{center}
\centering
    \begin{subfigure}[b]{0.4\textwidth}
        \centering
        \resizebox{\linewidth}{!}{
            \input{figure_warehouse_single}
        }
        \caption{Single-block layout.}
        \label{fig:warehouse_single}
    \end{subfigure}
    \begin{subfigure}[b]{0.4\textwidth}
        \centering
        \resizebox{\linewidth}{!}{
            \input{figure_warehouse_two}
        }
        \caption{Two-block layout.}
        \label{fig:warehouse_two}
    \end{subfigure}
\caption{
    Rectangular warehouse examples
    where the black squares represent the
    locations of items to be picked
    and the gray circle shows the depot.
} 
\label{fig:warehouse}
\end{center}
\vskip -0.2in
\end{figure}

%% file: figure_warehouse_single.tex
\begin{tikzpicture}[shorten >=1pt,draw=black!50]

    \pgfmathsetmacro{\x}{1}
    \pgfmathsetmacro{\y}{1.75}
    \pgfmathsetmacro{\m}{3}
    \pgfmathsetmacro{\n}{10} 
    \pgfmathsetmacro{\nodesize}{30}

    \foreach \j in {1,...,\m}{

        \draw[black, thin]
        (3*\j*\x-0.5*\x, 0.5*\x) --
        (3*\j*\x-0.5*\x, \n*\x + 0.5*\x) --
        (3*\j*\x-1.5*\x, \n*\x + 0.5*\x) --
        (3*\j*\x-1.5*\x, 0.5*\x) --
        cycle;
    
        \draw[black, thin]
        (3*\j*\x+0.5*\x, 0.5*\x) --
        (3*\j*\x+0.5*\x, \n*\x + 0.5*\x) --
        (3*\j*\x+1.5*\x, \n*\x + 0.5*\x) --
        (3*\j*\x+1.5*\x, 0.5*\x) --
        cycle;

        \node[rotate=90] (aisle-\j) at (3*\j*\x, 0.5*\n*\x+0.5*\x) {\Large aisle \j};

        \foreach \i in {2,...,\n}{
            \draw[black, thin] (3*\j*\x-0.5*\x, \i*\x-0.5*\x) -- (3*\j*\x-1.5*\x, \i*\x-0.5*\x);
            \draw[black, thin] (3*\j*\x+0.5*\x, \i*\x-0.5*\x) -- (3*\j*\x+1.5*\x, \i*\x-0.5*\x);

            \node (l-\j-\i) at (3*\j*\x-\x, \i*\x) { };
            \node (r-\j-\i) at (3*\j*\x+\x, \i*\x) { };
        }
    }

    \node (cross_front) at (1.5*\m*\x+1.5*\x, 0) {\Large front cross-aisle};
    \node (cross_back) at (1.5*\m*\x+1.5*\x, \n*\x+\x) {\Large back cross-aisle};

    \node[circle,draw=black,fill=gray!30,minimum size=\nodesize pt] (depot) at (3*\x, 0) {};
    \node[rectangle,fill=black,minimum size=\nodesize pt] (item1) at (l-1-3) {};
    \node[rectangle,fill=black,minimum size=\nodesize pt] (item2) at (r-2-4) {};
    \node[rectangle,fill=black,minimum size=\nodesize pt] (item3) at (l-2-8) {};
    \node[rectangle,fill=black,minimum size=\nodesize pt] (item4) at (r-3-2) {};
    \node[rectangle,fill=black,minimum size=\nodesize pt] (item5) at (r-3-9) {};

    \node (left) at (\x, -0.75\x) { };
    \node (right) at (3*\m*\x+2*\x, \n*\x+1.5*\x) { };

\end{tikzpicture}

%% file: figure_warehouse_two.tex
\begin{tikzpicture}[shorten >=1pt,draw=black!50]

    \pgfmathsetmacro{\x}{1}
    \pgfmathsetmacro{\y}{1.75}
    \pgfmathsetmacro{\m}{3}
    \pgfmathsetmacro{\n}{11} 
    \pgfmathsetmacro{\nodesize}{30}

    \foreach \j in {1,...,\m}{

        \draw[black, thin]
        (3*\j*\x-0.5*\x, 0.5*\x) --
        (3*\j*\x-0.5*\x, 0.5*\n*\x) --
        (3*\j*\x-1.5*\x, 0.5*\n*\x) --
        (3*\j*\x-1.5*\x, 0.5*\x) --
        cycle;

        \draw[black, thin]
        (3*\j*\x-0.5*\x, 0.5*\n*\x + \x) --
        (3*\j*\x-0.5*\x, \n*\x + 0.5*\x) --
        (3*\j*\x-1.5*\x, \n*\x + 0.5*\x) --
        (3*\j*\x-1.5*\x, 0.5*\n*\x + \x) --
        cycle;
    
        \draw[black, thin]
        (3*\j*\x+0.5*\x, 0.5*\x) --
        (3*\j*\x+0.5*\x, 0.5*\n*\x) --
        (3*\j*\x+1.5*\x, 0.5*\n*\x) --
        (3*\j*\x+1.5*\x, 0.5*\x) --
        cycle;

        \draw[black, thin]
        (3*\j*\x+0.5*\x, 0.5*\n*\x + \x) --
        (3*\j*\x+0.5*\x, \n*\x + 0.5*\x) --
        (3*\j*\x+1.5*\x, \n*\x + 0.5*\x) --
        (3*\j*\x+1.5*\x, 0.5*\n*\x + \x) --
        cycle;

        \node[rotate=90]
        (aisle-\j) at (3*\j*\x, 0.28*\n*\x+0*\x) {\Large lower subaisle \j};

        \node[rotate=90]
        (aisle-\j) at (3*\j*\x, 0.72*\n*\x+\x) {\Large upper subaisle \j};

        \foreach \i in {2,...,\n}{
            \draw[black, thin] (3*\j*\x-0.5*\x, \i*\x-0.5*\x) -- (3*\j*\x-1.5*\x, \i*\x-0.5*\x);
            \draw[black, thin] (3*\j*\x+0.5*\x, \i*\x-0.5*\x) -- (3*\j*\x+1.5*\x, \i*\x-0.5*\x);

            \node (l-\j-\i) at (3*\j*\x-\x, \i*\x) { };
            \node (r-\j-\i) at (3*\j*\x+\x, \i*\x) { };
        }
    }

    \node (cross_front) at (1.5*\m*\x+1.5*\x, 0) {\Large front cross-aisle};
    \node (cross_back) at (1.5*\m*\x+1.5*\x, 0.5*\n*\x+0.5*\x) {\Large middle cross-aisle};
    \node (cross_back) at (1.5*\m*\x+1.5*\x, \n*\x+\x) {\Large back cross-aisle};

    \node[circle,draw=black,fill=gray!30,minimum size=\nodesize pt] (depot) at (3*\x, 0) {};
    \node[rectangle,fill=black,minimum size=\nodesize pt] (item1) at (l-1-3) {};
    \node[rectangle,fill=black,minimum size=\nodesize pt] (item2) at (r-2-4) {};
    \node[rectangle,fill=black,minimum size=\nodesize pt] (item3) at (l-2-9) {};
    \node[rectangle,fill=black,minimum size=\nodesize pt] (item4) at (r-3-2) {};
    \node[rectangle,fill=black,minimum size=\nodesize pt] (item5) at (r-3-10) {};

    \node (left) at (\x, -0.75\x) { };
    \node (right) at (3*\m*\x+2*\x, \n*\x+1.5*\x) { };

\end{tikzpicture}

%% file: figure_vertical.tex
\begin{figure}[ht]
\centering

\resizebox{0.4\textwidth}{!}{

\begin{tikzpicture}[shorten >=1pt,->,draw=black!50, node distance=\layersep]

\tikzset{minimum size=32pt}

\node[shape=circle,draw=black] (a1) at (0, 5) { };
\node at (a1) {$a_j$};
\node[shape=circle,draw=black] (d1) at (0, 3.5) { };
\node[shape=circle,draw=black] (c1) at (0, 1.5) { };
\node[shape=circle,draw=black] (b1) at (0, 0) {$b_j$}; 
\draw[-, thick, draw=black] (a1) -- (d1) -- (c1) -- (b1);

\node (i) at (0, -1) {(i)};

\node[shape=circle,draw=black] (a2) at (1.5, 5) { };
\node at (a2) {$a_j$};
\node[shape=circle,draw=black] (d2) at (1.5, 3.5) { };
\node[shape=circle,draw=black] (c2) at (1.5, 1.5) { };
\node[shape=circle,draw=black] (b2) at (1.5, 0) {$b_j$}; 
\draw[-, thick, draw=black, double, double distance between line centers=8pt] (a2) -- (d2) -- (c2);

\node (ii) at (1.5, -1) {(ii)};

\node[shape=circle,draw=black] (a3) at (3, 5) { };
\node at (a3) {$a_j$};
\node[shape=circle,draw=black] (d3) at (3, 3.5) { };
\node[shape=circle,draw=black] (c3) at (3, 1.5) { };
\node[shape=circle,draw=black] (b3) at (3, 0) {$b_j$}; 
\draw[-, thick, draw=black, double, double distance between line centers=8pt] (b3) -- (c3) -- (d3);

\node (iii) at (3, -1) {(iii)};

\node[shape=circle,draw=black] (a4) at (4.5, 5) { };
\node at (a4) {$a_j$};
\node[shape=circle,draw=black] (d4) at (4.5, 3.5) { };
\node[shape=circle,draw=black] (c4) at (4.5, 1.5) { };
\node[shape=circle,draw=black] (b4) at (4.5, 0) {$b_j$}; 
\draw[-, thick, draw=black, double, double distance between line centers=8pt] (a4) -- (d4);
\draw[-, thick, draw=black, double, double distance between line centers=8pt] (b4) -- (c4);

\node (iv) at (4.5, -1) {(iv)};

\node[shape=circle,draw=black] (a5) at (6, 5) { };
\node at (a5) {$a_j$};
\node[shape=circle,draw=black] (d5) at (6, 3.5) { };
\node[shape=circle,draw=black] (c5) at (6, 1.5) { };
\node[shape=circle,draw=black] (b5) at (6, 0) {$b_j$}; 
\draw[-, thick, draw=black, double, double distance between line centers=8pt] (a5) -- (d5) -- (c5) -- (b5);

\node (v) at (6, -1) {(v)};

\node[shape=circle,draw=black] (a6) at (7.5, 5) { };
\node at (a6) {$a_j$};
\node[shape=circle,draw=black] (d6) at (7.5, 3.5) { };
\node[shape=circle,draw=black] (c6) at (7.5, 1.5) { };
\node[shape=circle,draw=black] (b6) at (7.5, 0) {$b_j$}; 

\node (vi) at (7.5, -1) {(vi)};


\end{tikzpicture}

}

\caption{
    Valid vertical edge configurations
    within aisle $j$.
    Unlabeled nodes represent item
    vertices within the aisle.
} 
\label{fig:vertical}
\end{figure}
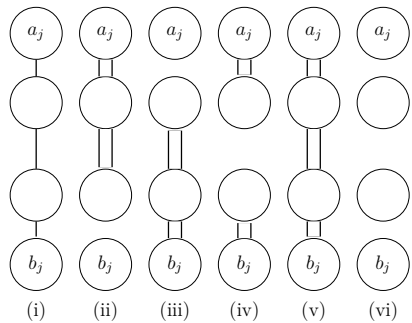

%% file: 3_single.tex

\section{Single-Block Rectangular Warehouses}
\label{single}

In this section,
we will briefly review
the single-block dynamic programming
algorithm of \citet{ratliff1983order}
and then propose a modification that preserves
the same state and action structure while reducing
the number of stages and transitions required.


The existing algorithm works from left to right,
applying all valid vertical edge configurations
within each aisle before applying
horizontal configurations between aisles.
The set of valid vertical edge configurations
are introduced in Section \ref{background}.
Between aisles $j$ and $j+1$,
there are five possible horizontal edge configurations,
\begin{equation*}
    11, 20, 02, 22, 00
\end{equation*}
where the first and second numbers indicate the number
of edges between the top and bottom cross-aisle vertices,
respectively (Figure \ref{fig:horizontal_single}).
The configuration $00$ is only valid if there are no items
to the right of the aisle $j$,
and is not valid elsewhere.

\input{figure_horizontal_single}

Two PTSs are equivalent if any subgraph that completes one
also completes the other.
An equivalence class is a set of equivalent PTSs
and is defined by the degree parity of the
rightmost cross-aisle vertices and the number of
connected components in each subgraph.
The states of the single-block algorithm
correspond to the seven PTS equivalence classes,
\begin{equation*}
    (U,U,1C), (0,E,1C), (E,0,1C), (E,E,1C), (E,E,2C), (0,0,0C), (0,0,1C)
\end{equation*}
where, for each class,
the first two components denote the
degree parity of the vertices $a_j$ and $b_j$, respectively
(zero $0$, odd $U$, or even $E$), and
the final component denotes connectivity 
($0C$, $1C$ or $2C$).


For a rectangular single-block warehouse, 
empty aisles may be ignored
and neither $2pass$ nor $none$ configurations are required
\cite{ratliff1983order, revenant2025note}.
Consequently,
only four valid vertical configurations need to be considered.
Since this reduction is already established in the literature,
we compare our proposed modified algorithm with a variant
of the original algorithm that omits these unnecessary actions.


\subsection{Original Single-Block Algorithm}


We briefly summarize the single-block dynamic programming
algorithm of \citet{ratliff1983order}.
A full description,
including the two state transition tables,
can be found in the original paper.

The algorithm begins at stage $L_1^-$
in the initial state $(0,0,0C)$ with cost zero.
In each aisle $j$,
a vertical stage is first performed with all valid
vertical configurations applied to the PTSs stored in $L_j^-$,
producing a set of $L_j^+$ PTSs.
From this set,
the minimal PTS is retained for each reachable state.
Next,
a horizontal stage is performed with all valid horizontal
configurations applied to the PTSs in $L_j^+$,
yielding the $L_{j+1}^-$ set,
again retaining only the minimal PTS for each state.

Each transition in the algorithm is of the form
\begin{equation*}
    (\text{state}, \text{ action}) \rightarrow \text{next state}
\end{equation*}
where an action corresponds to either a vertical
or horizontal configuration,
and the cost of the transition is equal to the total length
of the edges introduced by that configuration.
At every stage, the algorithm stores, for each state,
the action used to reach it,
the previous state,
and the accumulated cost.

After processing all $m$ aisles,
the optimal tour is obtained as the minimum PTS
in state $(0,0,1C)$ at stage $L_{m+1}^-$.
An example of an optimal tour is shown in
Figure \ref{fig:single_tour}.

\input{figure_single_tour}


\subsection{Modified Single-Block Algorithm}


We adapt the existing algorithm for rectangular warehouses
and demonstrate that vertical and horizontal stages
can be combined into a single stage.
\citet{dunn2025deterministic} proved that,
given the horizontal edges incident to an aisle,
the vertical edges required for a minimal tour
are uniquely determined.
In particular,
once the horizontal degrees of the upper and lower
cross-aisle vertices are known,
the minimal vertical configuration is 
$1pass$ if both degrees are odd;
$top$ if the back vertex is even and the front is not connected;
$bottom$ if the front vertex is even and the back is not connected;
and $gap$ if both vertices are even.

The implication for the existing algorithm is that,
given a minimal PTS in $L_j^-$,
the application of a horizontal configuration uniquely determines
the horizontal degree parities for aisle $j$ and therefore dictates
the necessary vertical configuration.
As a result,
the alternating vertical and horizontal stages of the
original algorithm can be collapsed into a single stage.
The modified algorithm proceeds by directly constructing
the minimal set of $L_{j+1}^-$ PTSs by applying all valid
horizontal configurations to the states in $L_j^-$,
together with the uniquely determined vertical configuration
implied by each horizontal choice.

Under this formulation, each transition is of the form
\begin{equation*}
    (\text{state, horizontal action})
    \rightarrow
    (\text{next state, vertical action)}
\end{equation*}
and the transition cost is equal to the total length of
the horizontal edges together with the required vertical edges.
Consequently, only a single state transition table is required
to define the algorithm,
incorporating both horizontal transitions
and associated vertical actions,
as shown in Table \ref{tab:single}.

\input{table_single}


\subsection{Single-Block Complexity Analysis}

\input{figure_stages_combined}

We compare the two algorithms in terms
of the total state-action evaluations performed.
For both algorithms,
the initial and terminal stages are restricted
to states $(0,0,0C)$ and $(0,0,1C)$, respectively,
reducing the number of valid transitions at
the boundaries.
Since these only occur once,
they do not affect the asymptotic complexity.
Accordingly,
our analysis focuses on the intermediate stages,
where the remaining five states are valid.

\paragraph{Original Algorithm}
At each vertical stage,
each of the five intermediate states has four valid actions,
resulting in $20$ state-action evaluations.
At each horizontal stage,
the number of valid actions varies by state with
a total of $9$ state-action evaluations
across all intermediate states.
Let $N$ denote the total number of state-action evaluations.
Since there are $m$ vertical stages and $m$
horizontal stages,
the total number of state-action evaluations
is approximately
\begin{equation*}
    N_{orig} \approx m(20+9) = 29m.
\end{equation*}

\paragraph{Modified Algorithm}
In the modified formulation,
each stage corresponds to a single combined transition
in which the vertical and horizontal actions that
previously occurred in separate stages are executed jointly.
At each stage,
the five intermediate states admit between two and four
valid horizontal actions,
resulting in a total of $16$ state–action evaluations per stage.
The total number of state–action evaluations is therefore approximately
\begin{equation*}
    N_{mod} \approx 16m.
\end{equation*}

\paragraph{Comparison}
A comparison of state space and transitions for
the single-block example is shown
in Figure \ref{fig:stages_combined}.
Unlike the original algorithm,
which requires $2m$ stages,
the modified algorithm requires only $m$ stages.
Both algorithms have linear time complexity in the number
of aisles,
\begin{equation*}
    \mathcal{O}(m),
\end{equation*}
since the state space and action sets are constant.
However, the modified algorithm reduces the number
of state-action evaluations by a factor of
\begin{equation*}
    \frac{N_{orig}}{N_{mod}} \approx \frac{29}{16} \approx 1.81,
\end{equation*}
corresponding to a reduction of approximately 45\%
in the number of evaluations.


\subsection{Single-Block Computational Experiments}

We evaluate the modified algorithm against the existing
single-block algorithm on randomly generated instances.
The experiments were implemented in Python on a laptop with
an Intel i7-13800H processor and 32 GB of RAM.
Benchmark instances were generated for all combinations of aisles
$m \in \{ 5, 10, 15, 20, 25, 30 \}$
and total pick positions
$|P| \in \{ 30, 45, 60, 75, 90 \}$
consistent with parameter settings
commonly used in the literature
\cite{scholz2016new, pansart2018exact, goeke2021modeling},
with 100 instances per combination.
The depot was randomly assigned to an aisle,
equally likely in the front or rear cross-aisle.
The number of pick locations per aisle was fixed at 90.
Action cost calculations were identical for both algorithms
and were excluded from runtime measurements


\paragraph{Speedup}

Table \ref{tab:single_speedup} reports the average speedup
(the ratio of original to modified runtime)
over the varying number of aisles and items.
Observed speedups range from 1.62 to 1.85,
with an overall average of 1.71.
The speedup remains fairly consistent
as the number of items changes,
which is expected since the computational complexity
of both algorithms is independent of the item count.
A mild decrease in speedup is observed
as the number of aisles increases,
which is consistent with
boundary stages and fixed overheads becoming
less dominant for larger instances.
Overall,
the empirical results closely match
the theoretical reduction in state–action evaluations.

\input{table_single_speedup}


\paragraph{Runtime}

Figure \ref{plt:single_runtime} plots the average runtime as
a function of the number of aisles.
As predicted by the complexity analysis,
runtimes appear to increase linearly for both methods,
with the modified algorithm consistently exhibiting
lower runtimes.

\input{plot_single_runtime}


%% file: figure_horizontal_single.tex
\begin{figure}[ht]
\centering
    \begin{minipage}{\textwidth}
        \centering
        \resizebox{0.8\textwidth}{!}
        {
        \begin{tikzpicture}
        [shorten >=1pt,->,draw=black!50, node distance=\layersep]

        \pgfmathsetmacro{\x}{1.75}
        \pgfmathsetmacro{\y}{3} 
        \pgfmathsetmacro{\nodesize}{30}


        \foreach \stage / \i in {0/(11),1/(20),2/(02),3/(22),4/(00)}{

            \node[shape=circle,draw=black, minimum size=\nodesize pt, fill=white]
            (a1-\stage) at (2*\stage*\x, \y) { };
            \node at (a1-\stage) {$a_j$};

            \node[shape=circle,draw=black, minimum size=\nodesize pt, fill=white]
            (a2-\stage) at (2*\stage*\x+1.1*\x, \y) { };
            \node at (a2-\stage) {$a_{j+1}$};

            \node[shape=circle,draw=black, minimum size=\nodesize pt, fill=white]
            (b1-\stage) at (2*\stage*\x, 0) { };
            \node at (b1-\stage) {$b_j$};

            \node[shape=circle,draw=black, minimum size=\nodesize pt, fill=white]
            (b2-\stage) at (2*\stage*\x+1.1*\x, 0) { };
            \node at (b2-\stage) {$b_{j+1}$};

            \node at (2*\stage*\x+0.5*\x, -0.5*\y) {\i};
            
        }


        \draw[-, thick, draw=black] 
        (a1-0) -- (a2-0);
        \draw[-, thick, draw=black] 
        (b1-0) -- (b2-0);

        \draw[-, thick, draw=black, double, double distance between line centers=8pt] 
        (a1-1) -- (a2-1);

        \draw[-, thick, draw=black, double, double distance between line centers=8pt] 
        (b1-2) -- (b2-2);

        \draw[-, thick, draw=black, double, double distance between line centers=8pt] 
        (a1-3) -- (a2-3);
        \draw[-, thick, draw=black, double, double distance between line centers=8pt] 
        (b1-3) -- (b2-3);

        
        \end{tikzpicture}
        }
        \caption{
            Valid horizontal edge configurations between
            aisle $j$ and $j+1$ of a single-block warehouse.
        } 
        \label{fig:horizontal_single}
    \end{minipage}
\end{figure}
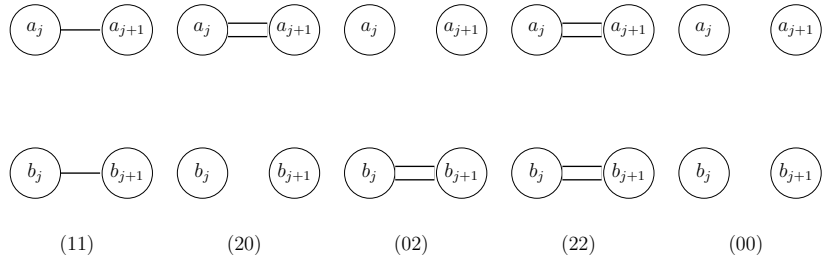

%% file: figure_single_tour.tex
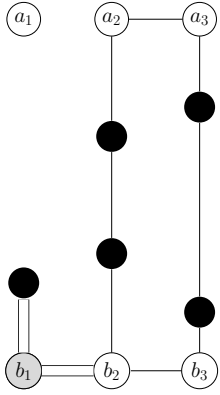
\begin{figure}[ht]
\centering
    \begin{minipage}{\textwidth}
        \centering
        \resizebox{0.3\textwidth}{!}
        {
        
        \begin{tikzpicture}[shorten >=1pt,draw=black!50]
        
            \pgfmathsetmacro{\x}{1}
            \pgfmathsetmacro{\y}{1.75}
            \pgfmathsetmacro{\m}{3}
            \pgfmathsetmacro{\n}{11}
            \pgfmathsetmacro{\nodesize}{30}
        
            \foreach \j in {1,...,\m}{
        
                \node[circle,draw=black,minimum size=\nodesize pt] (b-\j) at (3*\j*\x, 0) {\LARGE $b_{\j}$};
                \node[circle,draw=black,minimum size=\nodesize pt] (a-\j) at (3*\j*\x, \n*\x+\x) {\LARGE $a_{\j}$};
        
                \foreach \i in {2,...,\n}{
                    \node (\j-\i) at (3*\j*\x, \i*\x) { };
                }
            }
        
            \node[circle,draw=black,minimum size=\nodesize pt,fill=gray!30] (b-1) at (3*\x,0) {\LARGE $b_{1}$};
        
            \node[circle,fill=black,minimum size=\nodesize pt] (item1) at (1-3) {};
            \node[circle,fill=black,minimum size=\nodesize pt] (item2) at (2-4) {};
            \node[circle,fill=black,minimum size=\nodesize pt] (item3) at (2-8) {};
            \node[circle,fill=black,minimum size=\nodesize pt] (item4) at (3-2) {};
            \node[circle,fill=black,minimum size=\nodesize pt] (item5) at (3-9) {};
        
            \draw[-, thick, black, double, double distance between line centers=10pt]
            (b-2) -- (b-1) -- (item1);
        
            \draw[-, thick, black]
            (b-2) -- (item2) -- (item3) -- (a-2) -- (a-3)
            -- (item5) -- (item4) -- (b-3) -- (b-2);
        
            \node (left) at (\x, -0.75\x) { };
            \node (right) at (3*\m*\x+2*\x, \n*\x+1.5*\x) { };
        
        \end{tikzpicture}
        }
        \caption{
            Example of a minimal tour subgraph
            in a single-block rectangular warehouse.
            The solution has an optimal sequence of states            
            ($(0,0,0C)$,$(0,E,1C)$,$(0,E,1C)$,$(U,U,1C)$,
            $(U,U,1C)$,$(E,E,1C)$,$(0,0,1C)$)
            and actions
            $(bottom, 02, 1pass, 11, 1pass, 00)$.
        } 
        \label{fig:single_tour}
    \end{minipage}
\end{figure}

%% file: table_single.tex
\begin{table}[ht]
\centering
\caption {Single-block $L_{i+1}^-$ state transitions
and associated vertical edge configurations
that result from applying each horizontal configuration
to $L_i^-$ states.}
\label{tab:single}
\begin{center}
\begin{small}
\begin{sc}
    \resizebox{\textwidth}{!}{
    \begin{tabular}{| c | c  c  c  c  c  |}
    \hline 
    State & \multicolumn{5}{c|}{Horizontal Configurations} \\
    $L_i^-$ & $11$ & $20$ & $02$ & $22$ & $00^c$  \\
    \hline 
    $(U,U,1C)$ & $(U,U,1C)$ ($iv$) & $(E,0,1C)$ ($i$) & $(0,E,1C)$ ($i$) & $(E,E,1C)$ ($i$) & $(0,0,1C)$ ($i$) \\
    $(E,0,1C)$ & $(U,U,1C)$  ($i$) & $(E,0,1C)$ ($ii$) & $-$ & $(E,E,2C)$ ($iv$) & $(0,0,1C)$ ($ii$) \\
    $(0,E,1C)$ & $(U,U,1C)$  ($i$) & $-$ & $(0,E,1C)$ ($iii$) & $(E,E,2C)$ ($iv$) & $(0,0,1C)$ ($iii$) \\
    $(E,E,1C)$ & $(U,U,1C)$  ($i$) & $(E,0,1C)$ ($iv$) & $(0,E,1C)$ ($iv$) & $(E,E,1C)$ ($iv$) & $(0,0,1C)$ ($iv$) \\
    $(E,E,2C)$ & $(U,U,1C)$  ($i$) & $-$ & $-$ & $(E,E,2C)$ ($iv$) & $-$ \\
    $(0,0,0C)^a$ & $(U,U,1C)$  ($i$) & $(E,0,1C)$ ($ii$) & $(0,E,1C)$ ($iii$) & $(E,E,2C)$ ($iv$) & $(0,0,1C)$ ($iii$) \\
    $(0,0,1C)^b$ & $-$ & $-$ & $-$ & $-$ & $(0,0,1C)$ ($vi$) \\
    \hline
    \end{tabular}
    }
    \begin{tablenotes}
        \scriptsize
        \item $-$ Not valid or optimal.
        \item $^a$ Only possible in initial state.
        \item $^b$ Only possible in terminal state.
        \item $^c$ Only valid in final stage.
    \end{tablenotes}
\end{sc}
\end{small}
\end{center}
\end{table}

%% file: figure_stages_combined.tex
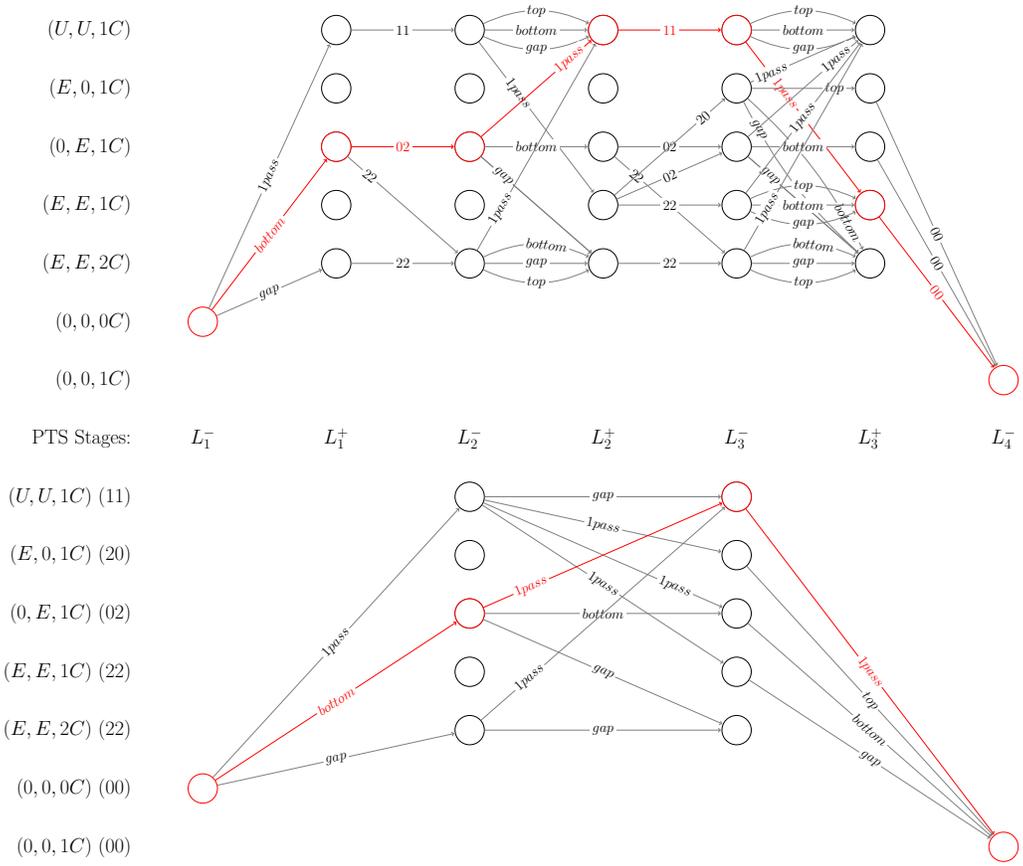
\begin{figure}[ht]
\centering
    \begin{minipage}{\textwidth}
        \centering
        \resizebox{\textwidth}{!}
        {
        \begin{tikzpicture}
        [shorten >=1pt,->,draw=black!50, node distance=\layersep]

        \tikzset{
          arrow label/.style={
            midway, sloped, fill=white, inner sep=2pt
          }
        }

        \tikzset{
          arrow left/.style={
            pos=0.2, sloped, fill=white, inner sep=2pt
          }
        }

        \tikzset{
          arrow right/.style={
            pos=0.8, sloped, fill=white, inner sep=2pt
          }
        }

        \pgfmathsetmacro{\x}{4}
        \pgfmathsetmacro{\y}{1.75} 
        \pgfmathsetmacro{\nodesize}{25}


        \node[anchor=east, font=\Large]
        (stages) at (-0.5*\x, -\y) {PTS Stages:};

        \node[anchor=east, font=\Large]
        (001C) at (-0.5*\x, 0) {$(0,0,1C)$};

        \node[anchor=east, font=\Large]
        (000C) at (-0.5*\x, \y) {$(0,0,0C)$};

        \node[anchor=east, font=\Large]
        (EE2C) at (-0.5*\x, 2*\y) {$(E,E,2C)$};

        \node[anchor=east, font=\Large]
        (EE1C) at (-0.5*\x, 3*\y) {$(E,E,1C)$};

        \node[anchor=east, font=\Large]
        (0E1C) at (-0.5*\x, 4*\y) {$(0,E,1C)$};

        \node[anchor=east, font=\Large]
        (E01C) at (-0.5*\x, 5*\y) {$(E,0,1C)$};

        \node[anchor=east, font=\Large]
        (UU1C) at (-0.5*\x, 6*\y) {$(U,U,1C)$};


        \node[shape=circle,draw=red, minimum size=\nodesize pt, fill=white]
        (000C-0) at (0, \y) {};

        \foreach \stage in {1,...,5}{
        


            \node[shape=circle,draw=black, minimum size=\nodesize pt, fill=white]
            (EE2C-\stage) at (\stage * \x, 2*\y) {};

            \node[shape=circle,draw=black, minimum size=\nodesize pt, fill=white]
            (EE1C-\stage) at (\stage * \x, 3*\y) {};

            \node[shape=circle,draw=black, minimum size=\nodesize pt, fill=white]
            (0E1C-\stage) at (\stage * \x, 4*\y) {};

            \node[shape=circle,draw=black, minimum size=\nodesize pt, fill=white]
            (E01C-\stage) at (\stage * \x, 5*\y) {};

            \node[shape=circle,draw=black, minimum size=\nodesize pt, fill=white]
            (UU1C-\stage) at (\stage * \x, 6*\y) {};
            
        }

        \node[shape=circle,draw=red, minimum size=\nodesize pt, fill=white]
        (001C-6) at (6*\x, 0) {};


        \node[font=\Large] (1-) at (0, -\y) {$L_1^-$};

        \draw[->] (000C-0) -- (UU1C-1) node[arrow label] {$1pass$};
        \draw[->, red] (000C-0) -- (0E1C-1) node[arrow label] {$bottom$};
        \draw[->] (000C-0) -- (EE2C-1) node[arrow label] {$gap$};


        \node[font=\Large] (1+) at (\x, -\y) {$L_1^+$};

        \node[shape=circle,draw=red, minimum size=\nodesize pt, fill=white]
        (0E1C-1) at (0E1C-1) {};

        \draw[->] (UU1C-1) -- (UU1C-2) node[arrow label] {$11$};
        
        \draw[->, red] (0E1C-1) -- (0E1C-2) node[arrow label] {$02$};
        \draw[->] (0E1C-1) -- (EE2C-2) node[arrow left] {$22$};
        
        
        
        \draw[->] (EE2C-1) -- (EE2C-2) node[arrow label] {$22$};


        \node[font=\Large] (2-) at (2*\x, -\y) {$L_2^-$};

        \node[shape=circle,draw=red, minimum size=\nodesize pt, fill=white]
        (0E1C-2) at (0E1C-2) {};

        \draw[->, bend left=25] (UU1C-2)
        to node[midway, sloped, fill=white, inner sep=2pt] {$top$} (UU1C-3);
        \draw[->] (UU1C-2)
        to node[midway, sloped, fill=white, inner sep=2pt] {$bottom$} (UU1C-3);
        \draw[->, bend right=25] (UU1C-2)
        to node[midway, sloped, fill=white, inner sep=2pt] {$gap$} (UU1C-3);
        \draw[->] (UU1C-2) -- (EE1C-3)
        node[pos=0.34, sloped, fill=white, inner sep=2pt] {$1pass$};


        \draw[->, red] (0E1C-2) -- (UU1C-3) node[arrow right] {$1pass$};
        \draw[->] (0E1C-2) -- (0E1C-3) node[arrow label] {$bottom$};
        \draw[->] (0E1C-2) -- (EE2C-3) node[arrow left] {$top$};
        \draw[->] (0E1C-2) -- (EE2C-3) node[arrow left] {$gap$};


        \draw[->] (EE2C-2) -- (EE2C-3) node[arrow label] {$gap$};
        \draw[->, bend right=25] (EE2C-2)
        to node[midway, sloped, fill=white, inner sep=2pt] {$top$} (EE2C-3);
        \draw[->, bend left=25] (EE2C-2)
        to node[pos=0.6, sloped, fill=white, inner sep=2pt] {$bottom$} (EE2C-3);
        \draw[->] (EE2C-2) -- (UU1C-3) node[arrow left] {$1pass$};


        \node[font=\Large] (2+) at (3*\x, -\y) {$L_2^+$};

        \node[shape=circle,draw=red, minimum size=\nodesize pt, fill=white]
        (UU1C-3) at (UU1C-3) {};

        \draw[->, red] (UU1C-3) -- (UU1C-4) node[arrow label] {$11$};
        
        \draw[->] (0E1C-3) -- (0E1C-4) node[arrow label] {$02$};
        \draw[->] (0E1C-3) -- (EE2C-4) node[arrow left] {$22$};
        
        
        \draw[->] (EE1C-3) -- (EE1C-4) node[arrow label] {$22$};
        \draw[->] (EE1C-3) -- (0E1C-4) node[arrow label] {$02$};
        \draw[->] (EE1C-3) -- (E01C-4) node[arrow right] {$20$};
        
        \draw[->] (EE2C-3) -- (EE2C-4) node[arrow label] {$22$};


        \node[font=\Large] (3-) at (4*\x, -\y) {$L_3^-$};

        \node[shape=circle,draw=red, minimum size=\nodesize pt, fill=white]
        (UU1C-4) at (UU1C-4) {};

        \draw[->, bend left=25] (UU1C-4)
        to node[midway, sloped, fill=white, inner sep=2pt] {$top$} (UU1C-5);
        \draw[->] (UU1C-4)
        to node[midway, sloped, fill=white, inner sep=2pt] {$bottom$} (UU1C-5);
        \draw[->, bend right=25] (UU1C-4)
        to node[midway, sloped, fill=white, inner sep=2pt] {$gap$} (UU1C-5);
        \draw[->, red] (UU1C-4) -- (EE1C-5)
        node[pos=0.34, sloped, fill=white, inner sep=2pt] {$1pass$};

        \draw[->] (E01C-4) -- (UU1C-5) node[arrow left] {$1pass$};
        \draw[->] (E01C-4) -- (E01C-5) node[arrow right] {$top$};
        \draw[->, bend left=10] (E01C-4)
        to node[pos=0.85, sloped, fill=white, inner sep=2pt] {$bottom$} (EE2C-5);
        \draw[->, bend right=10] (E01C-4)
        to node[pos=0.15, sloped, fill=white, inner sep=2pt] {$gap$} (EE2C-5);

        \draw[->] (0E1C-4) -- (UU1C-5) node[arrow right] {$1pass$};
        \draw[->] (0E1C-4) -- (0E1C-5) node[arrow label] {$bottom$};
        \draw[->] (0E1C-4) -- (EE2C-5) node[arrow left] {$top$};
        \draw[->] (0E1C-4) -- (EE2C-5) node[arrow left] {$gap$};

        \draw[->] (EE1C-4) -- (UU1C-5) node[arrow label] {$1pass$};
        \draw[->, bend left=25] (EE1C-4)
        to node[midway, sloped, fill=white, inner sep=2pt] {$top$} (EE1C-5);
        \draw[->] (EE1C-4)
        to node[midway, sloped, fill=white, inner sep=2pt] {$bottom$} (EE1C-5);
        \draw[->, bend right=25] (EE1C-4)
        to node[midway, sloped, fill=white, inner sep=2pt] {$gap$} (EE1C-5);

        \draw[->] (EE2C-4) -- (EE2C-5) node[arrow label] {$gap$};
        \draw[->, bend right=25] (EE2C-4)
        to node[midway, sloped, fill=white, inner sep=2pt] {$top$} (EE2C-5);
        \draw[->, bend left=25] (EE2C-4)
        to node[pos=0.6, sloped, fill=white, inner sep=2pt] {$bottom$} (EE2C-5);
        \draw[->] (EE2C-4) -- (UU1C-5) node[arrow left] {$1pass$};


        \node[font=\Large] (3+) at (5*\x, -\y) {$L_3^+$};

        \node[shape=circle,draw=red, minimum size=\nodesize pt, fill=white]
        (EE1C-5) at (EE1C-5) {};

        \draw[->] (E01C-5) -- (001C-6) node[arrow label] {$00$};
        \draw[->] (0E1C-5) -- (001C-6) node[arrow label] {$00$};
        \draw[->, red] (EE1C-5) -- (001C-6) node[arrow label] {$00$};


        \node[font=\Large] (4-) at (6*\x, -\y) {$L_4^-$};



        \node[anchor=east, font=\Large]
        (001Cn) at (-0.5*\x, -8*\y) {$(0,0,1C)$ $(00)$};

        \node[anchor=east, font=\Large]
        (000Cn) at (-0.5*\x, \y-8*\y) {$(0,0,0C)$ $(00)$};

        \node[anchor=east, font=\Large]
        (EE2Cn) at (-0.5*\x, 2*\y-8*\y) {$(E,E,2C)$ $(22)$};

        \node[anchor=east, font=\Large]
        (EE1Cn) at (-0.5*\x, 3*\y-8*\y) {$(E,E,1C)$ $(22)$};

        \node[anchor=east, font=\Large]
        (0E1Cn) at (-0.5*\x, 4*\y-8*\y) {$(0,E,1C)$ $(02)$};

        \node[anchor=east, font=\Large]
        (E01Cn) at (-0.5*\x, 5*\y-8*\y) {$(E,0,1C)$ $(20)$};

        \node[anchor=east, font=\Large]
        (UU1Cn) at (-0.5*\x, 6*\y-8*\y) {$(U,U,1C)$ $(11)$};


        \node[shape=circle,draw=red, minimum size=\nodesize pt, fill=white]
        (000Cn-0) at (0, \y-8*\y) {};

        \foreach \stage in {2, 4}{
        


            \node[shape=circle,draw=black, minimum size=\nodesize pt, fill=white]
            (EE2Cn-\stage) at (\stage * \x, 2*\y-8*\y) {};

            \node[shape=circle,draw=black, minimum size=\nodesize pt, fill=white]
            (EE1Cn-\stage) at (\stage * \x, 3*\y-8*\y) {};

            \node[shape=circle,draw=black, minimum size=\nodesize pt, fill=white]
            (0E1Cn-\stage) at (\stage * \x, 4*\y-8*\y) {};

            \node[shape=circle,draw=black, minimum size=\nodesize pt, fill=white]
            (E01Cn-\stage) at (\stage * \x, 5*\y-8*\y) {};

            \node[shape=circle,draw=black, minimum size=\nodesize pt, fill=white]
            (UU1Cn-\stage) at (\stage * \x, 6*\y-8*\y) {};
            
        }

        \node[shape=circle,draw=red, minimum size=\nodesize pt, fill=white]
        (001Cn-6) at (6*\x, 0-8*\y) {};


        \draw[->] (000Cn-0) -- (UU1Cn-2) node[arrow label] {$1pass$};
        \draw[->, red] (000Cn-0) -- (0E1Cn-2) node[arrow label] {$bottom$};
        \draw[->] (000Cn-0) -- (EE2Cn-2) node[arrow label] {$gap$};



        \node[shape=circle,draw=red, minimum size=\nodesize pt, fill=white]
        (0E1Cn-2) at (0E1Cn-2) {};

        \draw[->] (UU1Cn-2) -- (UU1Cn-4) node[arrow label] {$gap$};
        \draw[->] (UU1Cn-2) -- (E01Cn-4) node[arrow label] {$1pass$};
        \draw[->] (UU1Cn-2) -- (0E1Cn-4) node[arrow right] {$1pass$};
        \draw[->] (UU1Cn-2) -- (EE1Cn-4) node[arrow label] {$1pass$};


        \draw[->] (0E1Cn-2) -- (0E1Cn-4) node[arrow label] {$bottom$};
        \draw[->, red] (0E1Cn-2) -- (UU1Cn-4) node[arrow left] {$1pass$};
        \draw[->] (0E1Cn-2) -- (EE2Cn-4) node[arrow label] {$gap$};

        \draw[->] (EE2Cn-2) -- (UU1Cn-4) node[arrow left] {$1pass$};
        \draw[->] (EE2Cn-2) -- (EE2Cn-4) node[arrow label] {$gap$};



        \node[shape=circle,draw=red, minimum size=\nodesize pt, fill=white]
        (UU1Cn-4) at (UU1Cn-4) {};

        \draw[->, red] (UU1Cn-4) -- (001Cn-6) node[arrow label] {$1pass$};
        \draw[->] (E01Cn-4) -- (001Cn-6) node[arrow label] {$top$};
        \draw[->] (0E1Cn-4) -- (001Cn-6) node[arrow label] {$bottom$};
        \draw[->] (EE1Cn-4) -- (001Cn-6) node[arrow label] {$gap$};


        
        \end{tikzpicture}
        }
        \caption{State space of the original (top)
        and modified (bottom) algorithms for the
        warehouse in Figure \ref{fig:warehouse}.
        The minimal sequence of transitions and actions are
        highlighted in red.} 
        \label{fig:stages_combined}
    \end{minipage}
\end{figure}

%% file: table_single_speedup.tex
\begin{table}[t]
\centering
\caption{
    Average speedup of the modified algorithm ($t_{orig} / t_{mod}$)
    over rectangular single-block warehouse instances
    with varying number of aisles $m$
    and pick list sizes $|P|$.}
\label{tab:single_speedup}
\begin{filecontents*}{table_single_speedup.csv}
Items,5,10,15,20,25,30,Average
30,1.78,1.70,1.67,1.64,1.63,1.66,1.68
45,1.80,1.71,1.69,1.62,1.65,1.64,1.69
60,1.83,1.81,1.70,1.66,1.66,1.65,1.72
75,1.84,1.77,1.73,1.69,1.68,1.67,1.73
90,1.85,1.77,1.71,1.70,1.69,1.67,1.73
Avg.,1.82,1.75,1.70,1.66,1.66,1.66,1.71
\end{filecontents*}
\pgfplotstabletypeset[
 col sep=comma,
 string type,
 columns={Items,5,10,15,20,25,30,Average},
 columns/Items/.style={column name=$|P|$, column type=c},
 columns/5/.style={column name=5, column type=c},
 columns/10/.style={column name=10, column type=c},
 columns/15/.style={column name=15, column type=c},
 columns/20/.style={column name=20, column type=c},
 columns/25/.style={column name=25, column type=c},
 columns/30/.style={column name=30, column type=c},
 columns/Average/.style={column name=Avg., column type=c},
 every head row/.style={
  before row={
   \hline
   & \multicolumn{7}{c}{$m$} \\
   \cline{2-8}
  },
  after row=\hline
 },
 every last row/.style={after row=\hline}
]{table_single_speedup.csv}
\end{table}

%% file: plot_single_runtime.tex
\begin{figure}[t]
    \centering
    \begin{tikzpicture}
    \begin{axis}[
        width=0.95\linewidth,  
        height=0.5\linewidth, 
        xlabel={Number of aisles},
        ylabel={Runtime (\si{\milli\second})},
        legend pos=north west,
        grid=both,
        xtick={5,10,15,20,25,30},
        ytick={0.05,0.1,0.15},
        scaled y ticks=false,
        yticklabel style={/pgf/number format/fixed},
        mark size=2.5pt
    ]
    
    \addplot[black, thick, solid, mark=*] table[
        col sep=comma,
        x=n_aisles,
        y=original
    ]{plot_single_runtime.csv};
    \addlegendentry{Original}
    
    \addplot[black, thick, dashed, mark=square*, mark options={solid}] table[
        col sep=comma,
        x=n_aisles,
        y=modified
    ]{plot_single_runtime.csv};
    \addlegendentry{Modified}
    
    \end{axis}
    \end{tikzpicture}
\caption{
    Average runtime of the original
    and modified single-block algorithms.
}
\label{plt:single_runtime}
\end{figure}
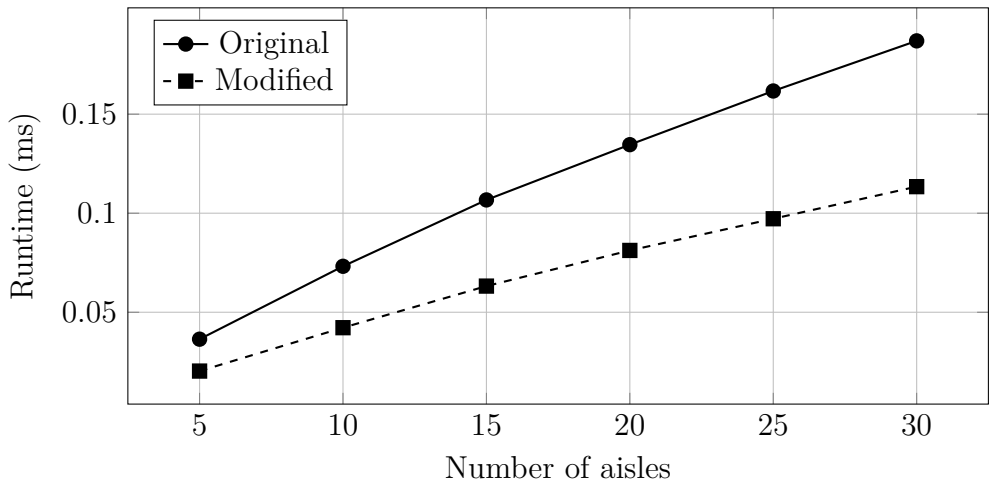

%% file: 4_two.tex

\section{Two-Block Rectangular Warehouses}
\label{two}

In this section,
we briefly review
the two-block algorithm of \citet{roodbergen2001routing}
and then propose two modifications that reduce
the number of required stages and transitions.
We begin by defining the actions and states
common to all three methods.


The vertical edge configurations
from the single-block algorithm apply again,
but are considered separately to the lower
and upper sub-aisles of each aisle.
For horizontal edges between aisles $j$ and $j+1$,
there are $14$ possible configurations:
\begin{equation*}
    110, 101, 011, 200, 020, 002, 211,
    121, 112, 220, 202, 022, 222, 000
\end{equation*}
where the first, second and third numbers denote the number
of edges between the top, middle
and bottom cross-aisle vertices,
respectively (Figure \ref{fig:horizontal_2block}).
The configuration $000$ is only valid if there are no items
to the right of aisle $j$.

\input{figure_horizontal_two}


The equivalence classes for a PTS in
a warehouse with three cross-aisles can be
characterized by
the degree parity of vertices $a_j$, $b_j$, and $c_j$;
the number of connected components;
and the distribution of $a_j$, $b_j$ and $c_j$
over the various components.
The last characteristic is only relevant when
all three vertices are even and the PTS
has two connected components;
for all other classes,
it can be ignored because the distribution is unique.
This yields $25$ equivalence classes for a two block warehouse
\begin{align*}
    &(E, 0, 0, 1C), (0, E, 0, 1C), (0, 0, E, 1C), (E, E, 0, 1C),\\
    &(E, 0, E, 1C), (0, E, E, 1C), (E, E, E, 1C), (U, U, 0, 1C),\\
    &(U, 0, U, 1C), (0, U, U, 1C), (E, U, U, 1C), (U, E, U, 1C),\\
    &(U, U, E, 1C), (E, E, 0, 2C), (E, 0, E, 2C), (0, E, E, 2C),\\
    &(E, E, E, 2C, a-bc), (E, E, E, 2C, b-ac), (E, E, E, 2C, c-ab),\\
    &(E, U, U, 2C), (U, E, U, 2C), (U, U, E, 2C), (E, E, E, 3C),\\
    &(0, 0, 0, 0C), (0, 0, 0, 1C).
\end{align*}



\subsection{Original Two-Block Algorithm}

We provide a brief overview of the two-block
dynamic programming algorithm of \citet{roodbergen2001routing}.
A full explanation,
including the three state transition tables,
can be found in the original paper.
As in the single-block scenario,
the algorithm moves sequentially from left to right
through the warehouse,
but now operates in two vertical stages,
one dedicated to the lower subaisle and one to the upper subaisle.

The algorithm begins at stage
$L_1^-$
in the initial state
$(0,0,0,0C)$
with cost zero.
In each aisle $j$,
PTSs are extended through a sequence of three stages.
First, vertical edges are added between the vertices
$b_j$ and $c_j$,
producing a set of candidate PTSs denoted $L_j^{+y}$.
Second, vertical edges are added between
the vertices $a_j$ and $b_j$,
which yields the set $L_j^{+x}$.
Finally, horizontal edges are added between the aisle vertices,
resulting in the set $L_{j+1}^-$.
At each stage,
only the minimum PTS is retained for each reachable state.
The optimal tour is the minimum PTS
in the terminal state $(0,0,0,1C)$.
An example of a minimal tour is shown in
Figure \ref{fig:two_tour}.

\input{figure_two_tour}


\subsection{Two-Block Algorithm Modification 1: No $2pass$}

Recent results show that double traversals are not required
in optimal tours for rectangular two-block warehouses with
more than one non-empty aisle \cite{dunn2025strict}.
Since this result was only recently 
established and is part of our own prior work,
we present the exclusion of $2pass$ vertical configurations
from the set of admissible actions in
the standard two-block algorithm as the
first modification.
This does not alter the state space
or the set of feasible tours,
but reduces the number of vertical actions
considered at each stage,
thus decreasing the number of state–action evaluations
required by the algorithm.

In contrast,
for the single-block case,
this restriction
is already well established in the literature
and was therefore incorporated directly
into the standard algorithm.


\subsection{Two-Block Algorithm Modification 2: Combined Stages}

We now extend the modification introduced
for single-block warehouses.
In a two-block layout,
applying horizontal edge configurations directly
to a $L_j^-$ PTS uniquely determines
the degree parity
and thus the necessary vertical
edge configurations for both the lower and upper subaisles.
All three stages can therefore be collapsed into one.
This reformulation preserves the original
state space and action feasibility,
but requires only a
single state transition table that simultaneously
encodes the horizontal transitions
and the associated vertical actions.
This formulation is shown in \ref{app1},
in which each transition takes the form
\begin{equation*}
    (\text{state, horizontal action})
    \rightarrow
    (\text{next state, lower action, upper action)}
\end{equation*}
where the transition cost is equal to the total length
of the horizontal edges,
together with the uniquely determined
vertical edges in both subaisles.
As in the single-block case,
only the minimal PTS is retained for each reachable state.


\subsection{Two-Block Complexity Analysis}

We compare the computational complexity of
the existing and modified two-block algorithms 
in terms of the total number of state–action 
evaluations performed. 
As in single-block algorithms,
the initial and terminal stages,
$(0,0,0,0C)$ and $(0,0,0,1C)$, respectively,
do not affect the asymptotic complexity.
The analysis therefore focuses on the intermediate stages,
where the remaining 23 states are admissible.


\paragraph{Original Algorithm}

For each of the three stages,
a table containing all state-action transitions
is provided in \citet{roodbergen2001routing}.
In the lower vertical stage
(between $b_j$ and $c_j$),
the 23 intermediate states admit a total of 130 valid actions.
In the upper vertical stage
(between $a_j$ and $b_j$),
there are again 130 state–action evaluations.
In the horizontal stage,
the 23 intermediate states admit 44 valid actions.
The total number of state–action evaluations
is therefore approximately
\begin{equation*}
    N_{orig} \approx m(130+130+44) = 304m.
\end{equation*}


\paragraph{Modification 1 (No $2pass$)}

With double traversals excluded,
the horizontal stage remains unchanged,
but the number of admissible vertical actions is reduced
in both vertical stages.
With the lower and upper vertical stages
each requiring 107 state–action evaluations,
the total number of state–action evaluations becomes
\begin{equation*}
    N_{mod1} \approx m(107+107+44) = 258m.
\end{equation*}


\paragraph{Modification 2 (Combined Stages)}

The two vertical stages are eliminated entirely
and all transitions are encoded within a single stage per aisle.
In this case,
each intermediate stage requires 231 state–action evaluations,
with the total number of evaluations approximately
\begin{equation*}
    N_{mod2} \approx 231m.
\end{equation*}


\paragraph{Comparison}

Although Modification 2 reduces 
the number of stages from $3m$ to $m$,
all three algorithms have linear time complexity
in the number of aisles,
$\mathcal{O}(m)$,
since the state space and action sets are constant.
However,
both modifications yield substantial
reductions in computational effort
arising first from
eliminating unnecessary vertical actions,
and then from collapsing multiple stages into
a single horizontal stage per aisle,
without altering the state space or the set of feasible tours.
This corresponds to a reduction by a factor of approximately 
$N_{orig} / N_{mod1} \approx 1.18$ and
$N_{orig} / N_{mod2} \approx 1.31$
relative to the original algorithm
for Modification 1 and Modification 2, respectively.


\subsection{Two-Block Computational Experiments}

This section evaluates the performance of the
proposed two-block algorithm modifications relative
to the original algorithm with double traversals.
Benchmark instances are the same as in the single-block experiments,
with the addition of a middle cross-aisle.
The cross-aisle has a width equal to the length of one
item location and is situated at the halfway point in
all aisles.


\paragraph{Speedup}

Table \ref{tab:two_speedup} reports the average speedup
of each modified algorithm relative to the original
for varying numbers of aisles.
The first modification yields a modest
but consistent improvement,
with an average speedup of 1.12,
which is close to the predicted reduction factor of 1.18.
The second modification achieves substantially larger gains,
with speedups ranging from 1.20 to 1.70 and an average of 1.31,
matching the theoretical reduction factor.
For both modifications,
the speedup decreases slightly as the number of aisles increases,
which is consistent with boundary stages and fixed overheads
becoming less significant for larger instances.

\input{table_two_speedup}


\paragraph{Runtime}

Figure \ref{plt:two_runtime} plots the average runtime
of all three algorithms.
As predicted by theoretical complexity,
all algorithms scale approximately linearly with the number of aisles.
Both modifications consistently outperform the original algorithm.
Modification 1 yields modest improvements,
while Modification 2
achieves substantial reductions in runtime,
with a stable speedup across all tested instance sizes.

\input{plot_two_runtime}

%% file: figure_horizontal_two.tex
\begin{figure}
\centering
\begin{center}

\def\v{1.25}
\def\h{1.5}
\def\d{5}

\resizebox{\textwidth}{!}{
    
\begin{tikzpicture} 

\tikzset{minimum size=25pt}

    \node[shape=circle,draw=black] (al1) at (0*\h, 2*\v) { };
    \node at (al1) {\small $a_j$}; 
    \node[shape=circle,draw=black] (ar1) at (1*\h, 2*\v) { };
    \node at (ar1) {\small $a_{j+1}$}; 
    \node[shape=circle,draw=black] (bl1) at (0*\h, 1*\v) { };
    \node at (bl1) {\small $b_j$}; 
    \node[shape=circle,draw=black] (br1) at (1*\h, 1*\v) { }; 
    \node at (br1) {\small $b_{j+1}$}; 
    \node[shape=circle,draw=black] (cl1) at (0*\h, 0*\v) { };
    \node at (cl1) {\small $c_j$}; 
    \node[shape=circle,draw=black] (cr1) at (1*\h, 0*\v) { }; 
    \node at (cr1) {\small $c_{j+1}$}; 
    \draw[-, thick] (al1) -- (ar1);
    \draw[-, thick] (cl1) -- (cr1);
    \draw[-, thick, double, double distance between line centers=\d pt] (bl1) -- (br1);
    \node (8) at (0.5*\h, -0.75*\v) {(121)};

    \node[shape=circle,draw=black] (al2) at (2*\h, 2*\v) { };
    \node at (al2) {\small $a_j$}; 
    \node[shape=circle,draw=black] (ar2) at (3*\h, 2*\v) { };
    \node at (ar2) {\small $a_{j+1}$}; 
    \node[shape=circle,draw=black] (bl2) at (2*\h, 1*\v) { };
    \node at (bl2) {\small $b_j$}; 
    \node[shape=circle,draw=black] (br2) at (3*\h, 1*\v) { }; 
    \node at (br2) {\small $b_{j+1}$}; 
    \node[shape=circle,draw=black] (cl2) at (2*\h, 0*\v) { };
    \node at (cl2) {\small $c_j$}; 
    \node[shape=circle,draw=black] (cr2) at (3*\h, 0*\v) { }; 
    \node at (cr2) {\small $c_{j+1}$}; 
    \draw[-, thick] (al2) -- (ar2);
    \draw[-, thick] (bl2) -- (br2);
    \draw[-, thick, double, double distance between line centers=\d pt] (cl2) -- (cr2);
    \node (9) at (2.5*\h, -0.75*\v) {(112)};

    \node[shape=circle,draw=black] (al3) at (4*\h, 2*\v) { };
    \node at (al3) {\small $a_j$}; 
    \node[shape=circle,draw=black] (ar3) at (5*\h, 2*\v) { };
    \node at (ar3) {\small $a_{j+1}$}; 
    \node[shape=circle,draw=black] (bl3) at (4*\h, 1*\v) { };
    \node at (bl3) {\small $b_j$}; 
    \node[shape=circle,draw=black] (br3) at (5*\h, 1*\v) { }; 
    \node at (br3) {\small $b_{j+1}$}; 
    \node[shape=circle,draw=black] (cl3) at (4*\h, 0*\v) { };
    \node at (cl3) {\small $c_j$}; 
    \node[shape=circle,draw=black] (cr3) at (5*\h, 0*\v) { }; 
    \node at (cr3) {\small $c_{j+1}$}; 
    \draw[-, thick, double, double distance between line centers=\d pt] (al3) -- (ar3);
    \draw[-, thick, double, double distance between line centers=\d pt] (bl3) -- (br3);
    \node (10) at (4.5*\h, -0.75*\v) {(220)};

    \node[shape=circle,draw=black] (al4) at (6*\h, 2*\v) { };
    \node at (al4) {\small $a_j$}; 
    \node[shape=circle,draw=black] (ar4) at (7*\h, 2*\v) { };
    \node at (ar4) {\small $a_{j+1}$}; 
    \node[shape=circle,draw=black] (bl4) at (6*\h, 1*\v) { };
    \node at (bl4) {\small $b_j$}; 
    \node[shape=circle,draw=black] (br4) at (7*\h, 1*\v) { }; 
    \node at (br4) {\small $b_{j+1}$}; 
    \node[shape=circle,draw=black] (cl4) at (6*\h, 0*\v) { };
    \node at (cl4) {\small $c_j$}; 
    \node[shape=circle,draw=black] (cr4) at (7*\h, 0*\v) { }; 
    \node at (cr4) {\small $c_{j+1}$}; 
    \draw[-, thick, double, double distance between line centers=\d pt] (al4) -- (ar4);
    \draw[-, thick, double, double distance between line centers=\d pt] (cl4) -- (cr4);
    \node (11) at (6.5*\h, -0.75*\v) {(202)};

    \node[shape=circle,draw=black] (al5) at (8*\h, 2*\v) { };
    \node at (al5) {\small $a_j$}; 
    \node[shape=circle,draw=black] (ar5) at (9*\h, 2*\v) { };
    \node at (ar5) {\small $a_{j+1}$}; 
    \node[shape=circle,draw=black] (bl5) at (8*\h, 1*\v) { };
    \node at (bl5) {\small $b_j$}; 
    \node[shape=circle,draw=black] (br5) at (9*\h, 1*\v) { }; 
    \node at (br5) {\small $b_{j+1}$}; 
    \node[shape=circle,draw=black] (cl5) at (8*\h, 0*\v) { };
    \node at (cl5) {\small $c_j$}; 
    \node[shape=circle,draw=black] (cr5) at (9*\h, 0*\v) { }; 
    \node at (cr5) {\small $c_{j+1}$}; 
    \draw[-, thick, double, double distance between line centers=\d pt] (bl5) -- (br5);
    \draw[-, thick, double, double distance between line centers=\d pt] (cl5) -- (cr5);
    \node (12) at (8.5*\h, -0.75*\v) {(022)};

    \node[shape=circle,draw=black] (al6) at (10*\h, 2*\v) { };
    \node at (al6) {\small $a_j$}; 
    \node[shape=circle,draw=black] (ar6) at (11*\h, 2*\v) { };
    \node at (ar6) {\small $a_{j+1}$}; 
    \node[shape=circle,draw=black] (bl6) at (10*\h, 1*\v) { };
    \node at (bl6) {\small $b_j$}; 
    \node[shape=circle,draw=black] (br6) at (11*\h, 1*\v) { }; 
    \node at (br6) {\small $b_{j+1}$}; 
    \node[shape=circle,draw=black] (cl6) at (10*\h, 0*\v) { };
    \node at (cl6) {\small $c_j$}; 
    \node[shape=circle,draw=black] (cr6) at (11*\h, 0*\v) { }; 
    \node at (cr6) {\small $c_{j+1}$}; 
    \draw[-, thick, double, double distance between line centers=\d pt] (al6) -- (ar6);
    \draw[-, thick, double, double distance between line centers=\d pt] (bl6) -- (br6);
    \draw[-, thick, double, double distance between line centers=\d pt] (cl6) -- (cr6);
    \node (13) at (10.5*\h, -0.75*\v) {(222)};

    \node[shape=circle,draw=black] (al7) at (12*\h, 2*\v) { };
    \node at (al7) {\small $a_j$}; 
    \node[shape=circle,draw=black] (ar7) at (13*\h, 2*\v) { };
    \node at (ar7) {\small $a_{j+1}$}; 
    \node[shape=circle,draw=black] (bl7) at (12*\h, 1*\v) { };
    \node at (bl7) {\small $b_j$}; 
    \node[shape=circle,draw=black] (br7) at (13*\h, 1*\v) { }; 
    \node at (br7) {\small $b_{j+1}$}; 
    \node[shape=circle,draw=black] (cl7) at (12*\h, 0*\v) { };
    \node at (cl7) {\small $c_j$}; 
    \node[shape=circle,draw=black] (cr7) at (13*\h, 0*\v) { }; 
    \node at (cr7) {\small $c_{j+1}$}; 
    \node (14) at (12.5*\h, -0.75*\v) {(000)};


    \node[shape=circle,draw=black] (al8) at (0*\h, 6*\v) { };
    \node at (al8) {\small $a_j$}; 
    \node[shape=circle,draw=black] (ar8) at (1*\h, 6*\v) { };
    \node at (ar8) {\small $a_{j+1}$}; 
    \node[shape=circle,draw=black] (bl8) at (0*\h, 5*\v) { };
    \node at (bl8) {\small $b_j$}; 
    \node[shape=circle,draw=black] (br8) at (1*\h, 5*\v) { }; 
    \node at (br8) {\small $b_{j+1}$}; 
    \node[shape=circle,draw=black] (cl8) at (0*\h, 4*\v) { };
    \node at (cl8) {\small $c_j$}; 
    \node[shape=circle,draw=black] (cr8) at (1*\h, 4*\v) { }; 
    \node at (cr8) {\small $c_{j+1}$}; 
    \draw[-, thick] (al8) -- (ar8);
    \draw[-, thick] (bl8) -- (br8);
    \node (1) at (0.5*\h, 3.25*\v) {(110)};

    \node[shape=circle,draw=black] (al9) at (2*\h, 6*\v) { };
    \node at (al9) {\small $a_j$}; 
    \node[shape=circle,draw=black] (ar9) at (3*\h, 6*\v) { };
    \node at (ar9) {\small $a_{j+1}$}; 
    \node[shape=circle,draw=black] (bl9) at (2*\h, 5*\v) { };
    \node at (bl9) {\small $b_j$}; 
    \node[shape=circle,draw=black] (br9) at (3*\h, 5*\v) { }; 
    \node at (br9) {\small $b_{j+1}$}; 
    \node[shape=circle,draw=black] (cl9) at (2*\h, 4*\v) { };
    \node at (cl9) {\small $c_j$}; 
    \node[shape=circle,draw=black] (cr9) at (3*\h, 4*\v) { }; 
    \node at (cr9) {\small $c_{j+1}$}; 
    \draw[-, thick] (al9) -- (ar9);
    \draw[-, thick] (cl9) -- (cr9);
    \node (2) at (2.5*\h, 3.25*\v) {(101)};

    \node[shape=circle,draw=black] (al10) at (4*\h, 6*\v) { };
    \node at (al10) {\small $a_j$}; 
    \node[shape=circle,draw=black] (ar10) at (5*\h, 6*\v) { };
    \node at (ar10) {\small $a_{j+1}$}; 
    \node[shape=circle,draw=black] (bl10) at (4*\h, 5*\v) { };
    \node at (bl10) {\small $b_j$}; 
    \node[shape=circle,draw=black] (br10) at (5*\h, 5*\v) { }; 
    \node at (br10) {\small $b_{j+1}$}; 
    \node[shape=circle,draw=black] (cl10) at (4*\h, 4*\v) { };
    \node at (cl10) {\small $c_j$}; 
    \node[shape=circle,draw=black] (cr10) at (5*\h, 4*\v) { }; 
    \node at (cr10) {\small $c_{j+1}$}; 
    \draw[-, thick] (bl10) -- (br10);
    \draw[-, thick] (cl10) -- (cr10);
    \node (3) at (4.5*\h, 3.25*\v) {(011)};

    \node[shape=circle,draw=black] (al11) at (6*\h, 6*\v) { };
    \node at (al11) {\small $a_j$}; 
    \node[shape=circle,draw=black] (ar11) at (7*\h, 6*\v) { };
    \node at (ar11) {\small $a_{j+1}$}; 
    \node[shape=circle,draw=black] (bl11) at (6*\h, 5*\v) { };
    \node at (bl11) {\small $b_j$}; 
    \node[shape=circle,draw=black] (br11) at (7*\h, 5*\v) { }; 
    \node at (br11) {\small $b_{j+1}$}; 
    \node[shape=circle,draw=black] (cl11) at (6*\h, 4*\v) { };
    \node at (cl11) {\small $c_j$}; 
    \node[shape=circle,draw=black] (cr11) at (7*\h, 4*\v) { }; 
    \node at (cr11) {\small $c_{j+1}$}; 
    \draw[-, thick, double, double distance between line centers=\d pt] (al11) -- (ar11);
    \node (4) at (6.5*\h, 3.25*\v) {(200)};

    \node[shape=circle,draw=black] (al12) at (8*\h, 6*\v) { };
    \node at (al12) {\small $a_j$}; 
    \node[shape=circle,draw=black] (ar12) at (9*\h, 6*\v) { };
    \node at (ar12) {\small $a_{j+1}$}; 
    \node[shape=circle,draw=black] (bl12) at (8*\h, 5*\v) { };
    \node at (bl12) {\small $b_j$}; 
    \node[shape=circle,draw=black] (br12) at (9*\h, 5*\v) { }; 
    \node at (br12) {\small $b_{j+1}$}; 
    \node[shape=circle,draw=black] (cl12) at (8*\h, 4*\v) { };
    \node at (cl12) {\small $c_j$}; 
    \node[shape=circle,draw=black] (cr12) at (9*\h, 4*\v) { }; 
    \node at (cr12) {\small $c_{j+1}$}; 
    \draw[-, thick, double, double distance between line centers=\d pt] (bl12) -- (br12);
    \node (5) at (8.5*\h, 3.25*\v) {(020)};

    \node[shape=circle,draw=black] (al13) at (10*\h, 6*\v) { };
    \node at (al13) {\small $a_j$}; 
    \node[shape=circle,draw=black] (ar13) at (11*\h, 6*\v) { };
    \node at (ar13) {\small $a_{j+1}$}; 
    \node[shape=circle,draw=black] (bl13) at (10*\h, 5*\v) { };
    \node at (bl13) {\small $b_j$}; 
    \node[shape=circle,draw=black] (br13) at (11*\h, 5*\v) { }; 
    \node at (br13) {\small $b_{j+1}$}; 
    \node[shape=circle,draw=black] (cl13) at (10*\h, 4*\v) { };
    \node at (cl13) {\small $c_j$}; 
    \node[shape=circle,draw=black] (cr13) at (11*\h, 4*\v) { }; 
    \node at (cr13) {\small $c_{j+1}$}; 
    \draw[-, thick, double, double distance between line centers=\d pt] (cl13) -- (cr13);
    \node (6) at (10.5*\h, 3.25*\v) {(002)};

    \node[shape=circle,draw=black] (al14) at (12*\h, 6*\v) { };
    \node at (al14) {\small $a_j$}; 
    \node[shape=circle,draw=black] (ar14) at (13*\h, 6*\v) { };
    \node at (ar14) {\small $a_{j+1}$}; 
    \node[shape=circle,draw=black] (bl14) at (12*\h, 5*\v) { };
    \node at (bl14) {\small $b_j$}; 
    \node[shape=circle,draw=black] (br14) at (13*\h, 5*\v) { }; 
    \node at (br14) {\small $b_{j+1}$}; 
    \node[shape=circle,draw=black] (cl14) at (12*\h, 4*\v) { };
    \node at (cl14) {\small $c_j$}; 
    \node[shape=circle,draw=black] (cr14) at (13*\h, 4*\v) { }; 
    \node at (cr14) {\small $c_{j+1}$}; 
    \draw[-, thick] (bl14) -- (br14);
    \draw[-, thick] (cl14) -- (cr14);
    \draw[-, thick, double, double distance between line centers=\d pt] (al14) -- (ar14);
    \node (7) at (12.5*\h, 3.25*\v) {(211)};

\end{tikzpicture}
}


\caption{Possible horizontal edge configurations.}
\label{fig:horizontal_2block}
\end{center}
\vskip -0.2in
\end{figure}
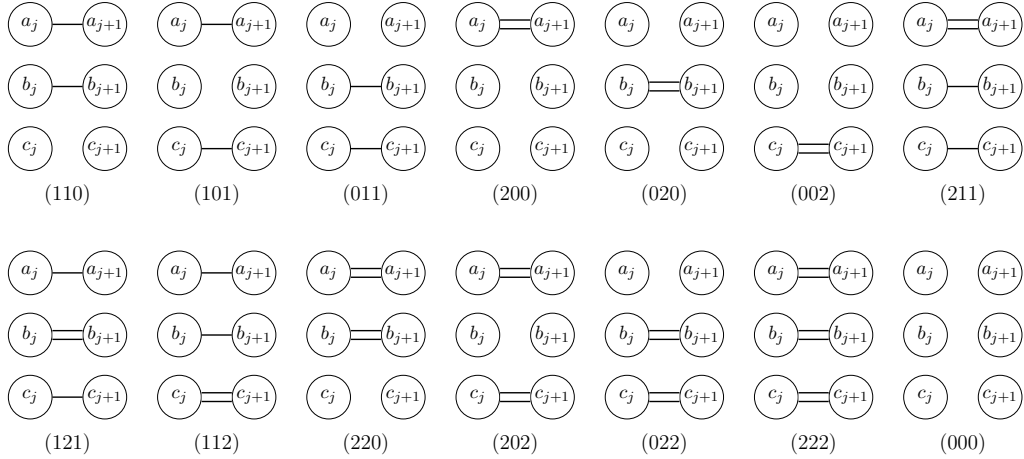

%% file: figure_two_tour.tex
\begin{figure}[ht]
\centering
    \begin{minipage}{\textwidth}
        \centering
        \resizebox{0.3\textwidth}{!}
        {
        
        \begin{tikzpicture}[shorten >=1pt,draw=black!50]
        
            \pgfmathsetmacro{\x}{1}
            \pgfmathsetmacro{\y}{1.75}
            \pgfmathsetmacro{\m}{3}
            \pgfmathsetmacro{\n}{12}
            \pgfmathsetmacro{\nodesize}{30}
        
            \foreach \j in {1,...,\m}{
        
                \node[circle,draw=black,minimum size=\nodesize pt] (c-\j) at (3*\j*\x, 0) {\LARGE $c_{\j}$};
                \node[circle,draw=black,minimum size=\nodesize pt] (b-\j) at (3*\j*\x, 0.5*\n*\x+0.5*\x) {\LARGE $b_{\j}$};
                \node[circle,draw=black,minimum size=\nodesize pt] (a-\j) at (3*\j*\x, \n*\x+\x) {\LARGE $a_{\j}$};
        
                \foreach \i in {2,...,\n}{
                    \node (\j-\i) at (3*\j*\x, \i*\x) { };
                }
            }
        
            \node[circle,draw=black,minimum size=\nodesize pt,fill=gray!30] (c-1) at (3*\x,0) {\LARGE $c_{1}$};
        
            \node[circle,fill=black,minimum size=\nodesize pt] (item1) at (1-3) {};
            \node[circle,fill=black,minimum size=\nodesize pt] (item2) at (2-4) {};
            \node[circle,fill=black,minimum size=\nodesize pt] (item3) at (2-10) {};
            \node[circle,fill=black,minimum size=\nodesize pt] (item4) at (3-2) {};
            \node[circle,fill=black,minimum size=\nodesize pt] (item5) at (3-11) {};
        
            \draw[-, thick, black, double, double distance between line centers=10pt]
            (b-2) -- (item2);
        
            \draw[-, thick, black]
            (c-1) -- (item1) -- (b-1) -- (b-2) -- (item3)
            -- (a-2) -- (a-3) -- (item5) -- (b-3) -- (item4)
            -- (c-3) -- (c-2) -- (c-1);
        
            \node (left) at (\x, -0.75\x) { };
            \node (right) at (3*\m*\x+2*\x, \n*\x+1.5*\x) { };
        
        \end{tikzpicture}
        }
        \caption{
            Example of a minimal tour subgraph
            in a two-block rectangular warehouse.
            The solution has an optimal sequence of states
            ($(0,0,0,0C)$,
            $(U,U,0,1C)$, $(U,U,0,1C)$, $(U,U,0,1C)$,
            $(U,U,0,1C)$, $(U,E,U,1C)$, $(U,0,U,1C)$,
            $(U,U,E,1C)$, $(E,E,E,1C)$, $(0,0,0,1C)$)
            and actions
            ($1pass$, $none$, $011$, $top$,
            $1pass$, $101$, $1pass$, $1pass$).
        } 
        \label{fig:two_tour}
    \end{minipage}
\end{figure}
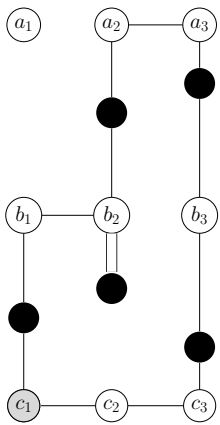

%% file: table_two_speedup.tex
\begin{table}[t]
\centering
\caption{
    Average speedup of each modified algorithm
    over rectangular two-block warehouse instances
    with varying number of aisles $m$.
}
\label{tab:two_speedup}

\begin{filecontents*}{table_two_speedup.csv}
algorithm,5,10,15,20,25,30,average
modified1,1.20,1.11,1.08,1.10,1.11,1.10,1.12
modified2,1.70,1.30,1.22,1.22,1.22,1.20,1.31
\end{filecontents*}
\pgfplotstabletypeset[
    col sep=comma,
    string type,
    columns={algorithm,5,10,15,20,25,30,average},
    columns/algorithm/.style={
        column name= ,
        column type=c,
        string replace={modified1}{$t_{orig} / t_{mod1}$},
        string replace={modified2}{$t_{orig} / t_{mod2}$}
    },
    columns/5/.style={column name=5,  column type=c},
    columns/10/.style={column name=10, column type=c},
    columns/15/.style={column name=15, column type=c},
    columns/20/.style={column name=20, column type=c},
    columns/25/.style={column name=25, column type=c},
    columns/30/.style={column name=30, column type=c},
    columns/average/.style={column name=Avg., column type=c},
    every head row/.style={
        before row={
            \hline
            & \multicolumn{7}{c}{$m$} \\
            \cline{2-8}
        },
        after row=\hline
    },
    every last row/.style={after row=\hline}
]{table_two_speedup.csv}

\end{table}

%% file: plot_two_runtime.tex
\begin{figure}[t]
    \centering
    \begin{tikzpicture}
    \begin{axis}[
        width=1.0\linewidth,  
        height=0.5\linewidth, 
        xlabel={Number of aisles},
        ylabel={Runtime (\si{\milli\second})},
        legend pos=north west,
        grid=both,
        xtick={5,10,15,20,25,30},
        scaled y ticks=false,
        yticklabel style={/pgf/number format/fixed},
        mark size=2.5pt
    ]
    
    \addplot[black, thick, mark=*] table[
        col sep=comma,
        x=n_aisles,
        y=original_2pass
    ]{plot_two_runtime.csv};
    \addlegendentry{Original}
    
    \addplot[black, thick, dotted, mark=x, mark options={solid}] table[
        col sep=comma,
        x=n_aisles,
        y=original_no2pass
    ]{plot_two_runtime.csv};
    \addlegendentry{Modification 1}

    \addplot[black, thick, dashed, mark=square*, mark options={solid}] table[
        col sep=comma,
        x=n_aisles,
        y=modified
    ]{plot_two_runtime.csv};
    \addlegendentry{Modification 2}
    
    \end{axis}
    \end{tikzpicture}
    \caption{
        Average runtime of the original
        and modified two-block algorithms.}
    \label{plt:two_runtime}
\end{figure}
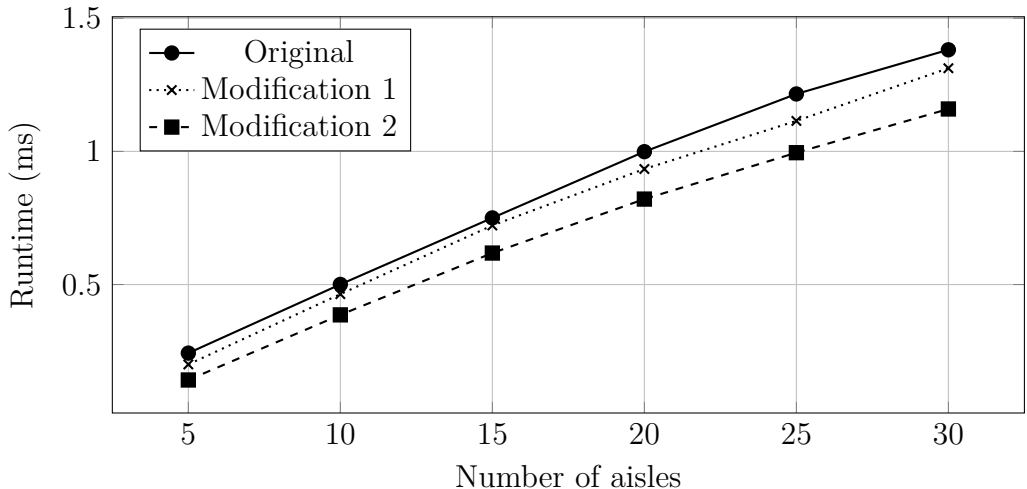

%% file: 5_discussion.tex

\section{Conclusion}
\label{discussion}


This paper leveraged recent structural
results on picker routing in
rectangular warehouses to propose
modifications to existing dynamic programming algorithms
that improve computational efficiency
without sacrificing optimality.
For both single-block and two-block layouts,
the proposed modifications preserve linear-time complexity
in the number of aisles while reducing the number
of state-action evaluations per stage.


Theoretical analysis predicts reductions
in computational effort of approximately 1.81
for the single-block algorithm and 1.18 and 1.31
for the two successive two-block modifications.
Extensive computational experiments confirm that
these theoretical improvements translate directly
into practical runtime gains,
with observed average speedups of 1.71, 1.12, and 1.31,
respectively.
This close agreement between theory and experiment
indicates that runtime is largely dominated by state–action
evaluations and confirms that the proposed modifications
deliver consistent and meaningful performance improvements


The proposed action and stage reductions
do not extend directly to larger warehouse layouts.
While traditional single-block and two-block algorithms
can be adapted to warehouses with additional cross-aisles,
the resulting state and action spaces grow rapidly
and quickly become computationally prohibitive
\cite{roodbergen2001routing}.
To address this,
the algorithm of \citet{pansart2018exact}
assigns horizontal configurations sequentially
to each cross-aisle segment,
requiring a transition for each section rather than
a single horizontal stage per aisle.
Consequently,
the stage-collapsing strategy considered here
does not apply directly.
Likewise,
the strict elimination of double traversals cannot be applied,
as such configurations have been shown to be necessary
in warehouses with more than three cross-aisles.
However,
double traversals are not required for connectivity
or within the upper and lower sub-aisles,
suggesting that restricted variants may still be feasible
\cite{dunn2025double, dunn2025strict}.
Exploring how these ideas can be selectively applied to
larger layouts while preserving optimality is a promising
direction for future research.

%% file: table_two_v2.tex
\footnotesize
\setlength{\tabcolsep}{2pt}
\renewcommand{\arraystretch}{1.05}

\begin{longtable}{@{}l*{5}{c}@{}}
\caption{Two-block state transitions and
    necessary vertical configurations.}
\label{tab:two-block-transition}
\\
\hline
\textbf{State} & \textbf{110} & \textbf{101} & \textbf{011} & \textbf{200} & \textbf{020} \\
\hline
\endfirsthead

\multicolumn{6}{c}%
{{\bfseries \tablename\ \thetable{} -- continued from previous page}} \\
\hline
\textbf{State} & \textbf{110} & \textbf{101} & \textbf{011} & \textbf{200} & \textbf{020} \\
\hline
\endhead

\hline \multicolumn{6}{|r|}{{Continued on next page}} \\ \hline
\endfoot

\begin{filecontents*}{table_two.csv}
State,110,101,011,200,020,002,211,121,112,220,202,022,222,000
{(E,0,0,1C)},{((U,U,0,1C),(i),(ii))},{((U,0,U,1C),(i),(i))},{--},{((E,0,0,1C),(ii),(vi))},{--},{--},{((E,U,U,2C),(iv),(i))},{((U,E,U,1C),(i),(i))},{((U,U,E,2C),(i),(iv))},{((E,E,0,2C),(iv),(ii))},{((E,0,E,2C),(ii),(iii))},{--},{((E,E,E,3C),(iv),(iv))},{((0,0,0,1C),(ii),(vi))}
{(0,E,0,1C)},{((U,U,0,1C),(i),(ii))},{((U,0,U,1C),(i),(i))},{((0,U,U,1C),(iii),(i))},{--},{((0,E,0,1C),(iii),(ii))},{--},{((E,U,U,2C),(iv),(i))},{((U,E,U,1C),(i),(i))},{((U,U,E,2C),(i),(iv))},{((E,E,0,2C),(iv),(ii))},{--},{((0,E,E,2C),(iii),(iv))},{((E,E,E,3C),(iv),(iv))},{((0,0,0,1C),(iii),(ii))}
{(0,0,E,1C)},{--},{((U,0,U,1C),(i),(i))},{((0,U,U,1C),(iii),(i))},{--},{--},{((0,0,E,1C),(vi),(iii))},{((E,U,U,2C),(iv),(i))},{((U,E,U,1C),(i),(i))},{((U,U,E,2C),(i),(iv))},{--},{((E,0,E,2C),(ii),(iii))},{((0,E,E,2C),(iii),(iv))},{((E,E,E,3C),(iv),(iv))},{((0,0,0,1C),(vi),(iii))}
{(E,E,0,1C)},{((U,U,0,1C),(i),(ii))},{((U,0,U,1C),(i),(i))},{((0,U,U,1C),(iv),(i))},{((E,0,0,1C),(iv),(ii))},{((0,E,0,1C),(iv),(ii))},{--},{((E,U,U,1C),(iv),(i))},{((U,E,U,1C),(i),(i))},{((U,U,E,2C),(i),(iv))},{((E,E,0,1C),(iv),(ii))},{((E,0,E,2C),(iv),(iv))},{((0,E,E,2C),(iv),(iv))},{((E,E,E,2C,c{--}ab),(iv),(iv))},{((0,0,0,1C),(iv),(ii))}
{(E,0,E,1C)},{((U,U,0,1C),(i),(iv))},{((U,0,U,1C),(i),(i))},{((0,U,U,1C),(iv),(i))},{((E,0,0,1C),(ii),(iii))},{--},{((0,0,E,1C),(ii),(iii))},{((E,U,U,1C),(iv),(i))},{((U,E,U,1C),(i),(i))},{((U,U,E,1C),(i),(iv))},{((E,E,0,2C),(iv),(iv))},{((E,0,E,1C),(ii),(iii))},{((0,E,E,2C),(iv),(iv))},{((E,E,E,2C,b{--}ac),(iv),(iv))},{((0,0,0,1C),(ii),(iii))}
{(0,E,E,1C)},{((U,U,0,1C),(i),(iv))},{((U,0,U,1C),(i),(i))},{((0,U,U,1C),(iii),(i))},{--},{((0,E,0,1C),(iii),(iv))},{((0,0,E,1C),(iii),(iv))},{((E,U,U,2C),(iv),(i))},{((U,E,U,1C),(i),(i))},{((U,U,E,1C),(i),(iv))},{((E,E,0,2C),(iv),(iv))},{((E,0,E,2C),(iv),(iv))},{((0,E,E,1C),(iii),(iv))},{((E,E,E,2C,a{--}bc),(iv),(iv))},{((0,0,0,1C),(iii),(iv))}
{(E,E,E,1C)},{((U,U,0,1C),(i),(iv))},{((U,0,U,1C),(i),(i))},{((0,U,U,1C),(iv),(i))},{((E,0,0,1C),(iv),(iv))},{((0,E,0,1C),(iv),(iv))},{((0,0,E,1C),(iv),(iv))},{((E,U,U,1C),(iv),(i))},{((U,E,U,1C),(i),(i))},{((U,U,E,1C),(i),(iv))},{((E,E,0,1C),(iv),(iv))},{((E,0,E,1C),(iv),(iv))},{((0,E,E,1C),(iv),(iv))},{((E,E,E,1C),(iv),(iv))},{((0,0,0,1C),(iv),(iv))}
{(U,U,0,1C)},{((U,U,0,1C),(iv),(ii))},{((U,0,U,1C),(iv),(i))},{((0,U,U,1C),(i),(i))},{((E,0,0,1C),(i),(ii))},{((0,E,0,1C),(i),(ii))},{--},{((E,U,U,1C),(i),(i))},{((U,E,U,1C),(iv),(i))},{((U,U,E,2C),(iv),(iv))},{((E,E,0,1C),(i),(ii))},{((E,0,E,2C),(i),(iv))},{((0,E,E,2C),(i),(iv))},{((E,E,E,2C,c{--}ab),(i),(iv))},{((0,0,0,1C),(i),(ii))}
{(U,0,U,1C)},{((U,U,0,1C),(iv),(i))},{((U,0,U,1C),(ii),(iii))},{((0,U,U,1C),(i),(iv))},{((E,0,0,1C),(i),(i))},{((0,E,0,1C),(i),(i))},{((0,0,E,1C),(i),(i))},{((E,U,U,1C),(i),(iv))},{((U,E,U,2C),(iv),(iv))},{((U,U,E,1C),(iv),(i))},{((E,E,0,1C),(i),(i))},{((E,0,E,1C),(i),(i))},{((0,E,E,1C),(i),(i))},{((E,E,E,1C),(i),(i))},{((0,0,0,1C),(i),(i))}
{(0,U,U,1C)},{((U,U,0,1C),(i),(i))},{((U,0,U,1C),(i),(iv))},{((0,U,U,1C),(iii),(iv))},{--},{((0,E,0,1C),(iii),(i))},{((0,0,E,1C),(iii),(i))},{((E,U,U,2C),(iv),(iv))},{((U,E,U,1C),(i),(iv))},{((U,U,E,1C),(i),(i))},{((E,E,0,2C),(iv),(i))},{((E,0,E,2C),(iv),(i))},{((0,E,E,1C),(iii),(i))},{((E,E,E,2C,a{--}bc),(iv),(i))},{((0,0,0,1C),(iii),(i))}
{(E,U,U,1C)},{((U,U,0,1C),(i),(i))},{((U,0,U,1C),(i),(iv))},{((0,U,U,1C),(iv),(iv))},{((E,0,0,1C),(iv),(i))},{((0,E,0,1C),(iv),(i))},{((0,0,E,1C),(iv),(i))},{((E,U,U,1C),(iv),(iv))},{((U,E,U,1C),(i),(iv))},{((U,U,E,1C),(i),(i))},{((E,E,0,1C),(iv),(i))},{((E,0,E,1C),(iv),(i))},{((0,E,E,1C),(iv),(i))},{((E,E,E,1C),(iv),(i))},{((0,0,0,1C),(iv),(i))}
{(U,E,U,1C)},{((U,U,0,1C),(iv),(i))},{((U,0,U,1C),(iv),(iv))},{((0,U,U,1C),(i),(iv))},{((E,0,0,1C),(i),(i))},{((0,E,0,1C),(i),(i))},{((0,0,E,1C),(i),(i))},{((E,U,U,1C),(i),(iv))},{((U,E,U,1C),(iv),(iv))},{((U,U,E,1C),(iv),(i))},{((E,E,0,1C),(i),(i))},{((E,0,E,1C),(i),(i))},{((0,E,E,1C),(i),(i))},{((E,E,E,1C),(i),(i))},{((0,0,0,1C),(i),(i))}
{(U,U,E,1C)},{((U,U,0,1C),(iv),(iv))},{((U,0,U,1C),(iv),(i))},{((0,U,U,1C),(i),(i))},{((E,0,0,1C),(i),(iv))},{((0,E,0,1C),(i),(iv))},{((0,0,E,1C),(i),(iv))},{((E,U,U,1C),(i),(i))},{((U,E,U,1C),(iv),(i))},{((U,U,E,1C),(iv),(iv))},{((E,E,0,1C),(i),(iv))},{((E,0,E,1C),(i),(iv))},{((0,E,E,1C),(i),(iv))},{((E,E,E,1C),(i),(iv))},{((0,0,0,1C),(i),(iv))}
{(E,E,0,2C)},{((U,U,0,1C),(i),(ii))},{((U,0,U,1C),(i),(i))},{--},{--},{--},{--},{((E,U,U,2C),(iv),(i))},{((U,E,U,1C),(i),(i))},{((U,U,E,2C),(i),(iv))},{((E,E,0,2C),(iv),(ii))},{--},{--},{((E,E,E,3C),(iv),(iv))},{--}
{(E,0,E,2C)},{--},{((U,0,U,1C),(i),(i))},{--},{--},{--},{--},{((E,U,U,2C),(iv),(i))},{((U,E,U,1C),(i),(i))},{((U,U,E,2C),(i),(iv))},{--},{((E,0,E,2C),(ii),(iii))},{--},{((E,E,E,3C),(iv),(iv))},{--}
{(0,E,E,2C)},{--},{((U,0,U,1C),(i),(i))},{((0,U,U,1C),(iii),(i))},{--},{--},{--},{((E,U,U,2C),(iv),(i))},{((U,E,U,1C),(i),(i))},{((U,U,E,2C),(i),(iv))},{--},{--},{((0,E,E,2C),(iii),(iv))},{((E,E,E,3C),(iv),(iv))},{--}
{(E,E,E,2C,a{--}bc)},{((U,U,0,1C),(i),(iv))},{((U,0,U,1C),(i),(i))},{--},{--},{--},{--},{((E,U,U,2C),(iv),(i))},{((U,E,U,1C),(i),(i))},{((U,U,E,1C),(i),(iv))},{((E,E,0,2C),(iv),(iv))},{((E,0,E,2C),(iv),(iv))},{--},{((E,E,E,2C,a{--}bc),(iv),(iv))},{--}
{(E,E,E,2C,b{--}ac)},{((U,U,0,1C),(i),(iv))},{((U,0,U,1C),(i),(i))},{((0,U,U,1C),(iv),(i))},{--},{--},{--},{((E,U,U,1C),(iv),(i))},{((U,E,U,1C),(i),(i))},{((U,U,E,1C),(i),(iv))},{((E,E,0,2C),(iv),(iv))},{--},{((0,E,E,2C),(iv),(iv))},{((E,E,E,2C,b{--}ac),(iv),(iv))},{--}
{(E,E,E,2C,c{--}ab)},{--},{((U,0,U,1C),(i),(i))},{((0,U,U,1C),(iv),(i))},{--},{--},{--},{((E,U,U,1C),(iv),(i))},{((U,E,U,1C),(i),(i))},{((U,U,E,2C),(i),(iv))},{--},{((E,0,E,2C),(iv),(iv))},{((0,E,E,2C),(iv),(iv))},{((E,E,E,2C,c{--}ab),(iv),(iv))},{--}
{(E,U,U,2C)},{((U,U,0,1C),(i),(i))},{((U,0,U,1C),(i),(iv))},{--},{--},{--},{--},{((E,U,U,2C),(iv),(iv))},{((U,E,U,1C),(i),(iv))},{((U,U,E,1C),(i),(i))},{((E,E,0,2C),(iv),(i))},{((E,0,E,2C),(iv),(i))},{--},{((E,E,E,2C,a{--}bc),(iv),(i))},{--}
{(U,E,U,2C)},{((U,U,0,1C),(iv),(i))},{--},{((0,U,U,1C),(i),(iv))},{((E,0,0,1C),(i),(i))},{((0,E,0,1C),(i),(i))},{((0,0,E,1C),(i),(i))},{((E,U,U,1C),(i),(iv))},{((U,E,U,2C),(iv),(iv))},{((U,U,E,1C),(iv),(i))},{((E,E,0,1C),(i),(i))},{((E,0,E,1C),(i),(i))},{((0,E,E,1C),(i),(i))},{((E,E,E,1C),(i),(i))},{((0,0,0,1C),(i),(i))}
{(U,U,E,2C)},{--},{((U,0,U,1C),(iv),(i))},{((0,U,U,1C),(i),(i))},{--},{--},{--},{((E,U,U,1C),(i),(i))},{((U,E,U,1C),(iv),(i))},{((U,U,E,2C),(iv),(iv))},{--},{((E,0,E,2C),(i),(iv))},{((0,E,E,2C),(i),(iv))},{((E,E,E,2C,c{--}ab),(i),(iv))},{--}
{(E,E,E,3C)},{--},{((U,0,U,1C),(i),(i))},{--},{--},{--},{--},{((E,U,U,2C),(iv),(i))},{((U,E,U,1C),(i),(i))},{((U,U,E,2C),(i),(iv))},{--},{--},{--},{((E,E,E,3C),(iv),(iv))},{--}
{(0,0,0,0C)},{((U,U,0,1C),(i),(ii))},{((U,0,U,1C),(i),(i))},{((0,U,U,1C),(iii),(i))},{((E,0,0,1C),(ii),(vi))},{((0,E,0,1C),(iii),(ii))},{((0,0,E,1C),(vi),(iii))},{((E,U,U,2C),(iv),(i))},{((U,E,U,1C),(i),(i))},{((U,U,E,2C),(i),(iv))},{((E,E,0,2C),(iv),(ii))},{((E,0,E,2C),(ii),(iii))},{((0,E,E,2C),(iii),(iv))},{((E,E,E,3C),(iv),(iv))},{((0,0,0,0C),(vi),(vi))}
{(0,0,0,1C)},{--},{--},{--},{--},{--},{--},{--},{--},{--},{--},{--},{--},{--},{--}
\end{filecontents*}
\csvreader[
  head to column names,
  late after line=\\,
]{table_two.csv}{}{%
\csvcoli & \csvcolii & \csvcoliii & \csvcoliv & \csvcolv & \csvcolvi
}

\hline
\endlastfoot

\end{longtable}

\vspace{0.5em}

\begin{longtable}{@{}l*{5}{c}@{}}
\caption*{(continued)}
\\
\hline
\textbf{State} & \textbf{002} & \textbf{211} & \textbf{121} & \textbf{112} & \textbf{220} \\
\hline
\endfirsthead

\multicolumn{6}{c}%
{{\bfseries Table \thetable{} (Two-block transition table) -- continued from previous page}} \\
\hline
\textbf{State} & \textbf{002} & \textbf{211} & \textbf{121} & \textbf{112} & \textbf{220} \\
\hline
\endhead

\hline \multicolumn{6}{|r|}{{Continued on next page}} \\ \hline
\endfoot

\csvreader[
  head to column names,
  late after line=\\,
]{table_two.csv}{}{%
\csvcoli & \csvcolvii & \csvcolviii & \csvcolix & \csvcolx & \csvcolxi
}

\hline
\endlastfoot

\end{longtable}

\clearpage

\begin{longtable}{@{}l*{4}{c}@{}}
\caption*{(continued)}
\\
\hline
\textbf{State} & \textbf{202} & \textbf{022} & \textbf{222} & \textbf{000} \\
\hline
\endfirsthead

\multicolumn{5}{c}%
{{\bfseries Table \thetable{} (Two-block transition table) -- continued from previous page}} \\
\hline
\textbf{State} & \textbf{202} & \textbf{022} & \textbf{222} & \textbf{000} \\
\hline
\endhead

\hline \multicolumn{5}{|r|}{{Continued on next page}} \\ \hline
\endfoot

\csvreader[
  head to column names,
  late after line=\\,
]{table_two.csv}{}{%
\csvcoli & \csvcolxii & \csvcolxiii & \csvcolxiv & \csvcolxv
}

\hline
\endlastfoot

\end{longtable}

%% file: bibliography.bib
@article{dunn2025deterministic,
  title={Deterministic Structure of Vertical Configurations in Minimal Picker Tours for Rectangular Warehouses},
  author={Dunn, George and Stojanovski, Elizabeth and Lamichhane, Bishnu and Charkhgard, Hadi and Eshragh, Ali},
  journal={arXiv preprint arXiv:2508.00365},
  year={2025}
}

@article{dunn2025double,
  title={Double Traversals in Optimal Picker Routes for Warehouses with Multiple Blocks},
  journal = {Operations Research Letters},
  volume = {65},
  pages = {107397},
  year = {2026},
  issn = {0167-6377},
  doi = {https://doi.org/10.1016/j.orl.2025.107397},
  author = {George Dunn and Hadi Charkhgard and Ali Eshragh and Elizabeth Stojanovski}
}

@article{dunn2025strict,
  title={Strict Elimination of Double Traversals in Outer Subaisles and Two-Block Rectangular Warehouses},
  author={Dunn, George and Stojanovski, Elizabeth and Lamichhane, Bishnu and Charkhgard, Hadi and Eshragh, Ali},
  journal={arXiv preprint arXiv:2512.08235},
  year={2025}
}

@article{goeke2021modeling,
  title={Modeling single-picker routing problems in classical and modern warehouses: INFORMS journal on computing meritorious paper awardee},
  author={Goeke, Dominik and Schneider, Michael},
  journal={INFORMS Journal on Computing},
  volume={33},
  number={2},
  pages={436--451},
  year={2021},
  publisher={INFORMS}
}

@article{pansart2018exact,
    title = {Exact algorithms for the order picking problem},
    journal = {Computers \& Operations Research},
    volume = {100},
    pages = {117-127},
    year = {2018},
    issn = {0305-0548},
    doi = {https://doi.org/10.1016/j.cor.2018.07.002},
    author = {Lucie Pansart and Nicolas Catusse and Hadrien Cambazard}
}

@article{ratliff1983order,
    title={Order-picking in a rectangular warehouse: a solvable case of the traveling salesman problem},
    author={Ratliff, H Donald and Rosenthal, Arnon S},
    journal={Operations research},
    volume={31},
    number={3},
    pages={507--521},
    year={1983},
    publisher={INFORMS}
}

@article{revenant2025note,
    title={A note about a transition of Ratliff and Rosenthal's order picking algorithm for rectangular warehouses},
    author={Revenant, Paul and Cambazard, Hadrien and Catusse, Nicolas},
    journal={Operations Research Letters},
    volume={62},
    pages={107325},
    year={2025},
    publisher={Elsevier},
    doi={https://doi.org/10.1016/j.orl.2025.107325}
}

@article{roodbergen2001routing,
    title = {Routing order pickers in a warehouse with a middle aisle},
    journal = {European Journal of Operational Research},
    volume = {133},
    number = {1},
    pages = {32-43},
    year = {2001},
    issn = {0377-2217},
    doi = {https://doi.org/10.1016/S0377-2217(00)00177-6},
    author = {Kees Jan Roodbergen and René {de Koster}}
}

@article{scholz2016new,
  title={A new mathematical programming formulation for the single-picker routing problem},
  author={Scholz, Andr{\'e} and Henn, Sebastian and Stuhlmann, Meike and W{\"a}scher, Gerhard},
  journal={European Journal of Operational Research},
  volume={253},
  number={1},
  pages={68--84},
  year={2016},
  publisher={Elsevier}
}
